\documentclass[11pt,reqno,a4paper]{article}
\DeclareMathAlphabet{\mathpzc}{OT1}{pzc}{m}{it}

\usepackage{geometry}
\usepackage{amsmath,amsthm,amstext,amscd,amssymb,euscript,url}
\usepackage{mathrsfs}
\usepackage{calrsfs}
\usepackage{epsfig}
\usepackage[inline]{enumitem}
\usepackage{color}

\newcommand{\locb}{\mathcal{L}}

\newcommand{\col}[1]{\textcolor[rgb]{0,0,0}{#1}}

\newcommand{\rll}{\rho_\ell}
\newcommand{\rrr}{\rho_r}

\newcommand{\cll}{c_\ell}
\newcommand{\crr}{c_r}

\renewcommand{\k}{\kappa}

\newcommand{\Z}{\mathbb Z}
\newcommand{\R}{\mathbb R}
\newcommand{\N}{\mathbb N}

\newcommand{\E}{\mathbb E}

\renewcommand{\phi}{\varphi}

\newcommand{\si}{\ensuremath{\sigma}}

\newcommand{\pee}{\ensuremath{\mathbb{P}}}
\newcommand{\gee}{\ensuremath{\mathcal{G}}}
\newcommand{\loc}{\mathcal{L}}

\def\1{{\mathchoice {\rm 1\mskip-4mu l} {\rm 1\mskip-4mu l}
{\rm 1\mskip-4.5mu l} {\rm 1\mskip-5mu l}}}

\newtheorem{theorem}{{\small T}{\scriptsize HEOREM}}[section]
\newtheorem{corollary}{{\bf{\small C}{\scriptsize OROLLARY}}}[section]
\newtheorem{proposition}{{\bf{\small P}{\scriptsize ROPOSITION}}}[section]
\newtheorem{lemma}{{\bf{\small L}{\scriptsize EMMA}}}[section]
\newtheorem{remark}{{\bf{\small R}{\scriptsize EMARK}}}[section]
\newtheorem{definition}{{\bf{\small D}{\scriptsize EFINITION}}}[section]

\renewenvironment{proof}[1]
{\noindent{{\bf{\small{ P}{\scriptsize ROOF}}}.}\hspace{0.1cm} #1} {$\;\qed$\newline}

\newcommand{\beq}{\begin{eqnarray}}
\newcommand{\eeq}{\end{eqnarray}}

\newcommand{\ba}{\begin{align*}}
\newcommand{\ea}{\end{align*}}

\newcommand{\be}{\begin{equation}}
\newcommand{\ee}{\end{equation}}

\newcommand{\bl}{\begin{lemma}}
\newcommand{\el}{\end{lemma}}

\newcommand{\br}{\begin{remark}}
\newcommand{\er}{\end{remark}}

\newcommand{\bt}{\begin{theorem}}
\newcommand{\et}{\end{theorem}}

\newcommand{\bd}{\begin{definition}}
\newcommand{\ed}{\end{definition}}

\newcommand{\bp}{\begin{proposition}}
\newcommand{\ep}{\end{proposition}}

\newcommand{\bc}{\begin{corollary}}
\newcommand{\ec}{\end{corollary}}

\newcommand{\bpr}{\begin{proof}}
\newcommand{\epr}{\end{proof}}

\newcommand{\bi}{\begin{itemize}}
\newcommand{\ei}{\end{itemize}}

\newcommand{\ben}{\begin{enumerate}}
\newcommand{\een}{\end{enumerate}}

\newcommand{\nn}{\nonumber}

\newcommand{\caA}{{\mathcal A}}

\newcommand{\caC}{{\mathscr C}}

\newcommand{\caE}{{\mathrsfs E}}

\newcommand{\caL}{{\mathcal L}}

\newcommand{\caS}{{\mathcal S}}

\newcommand{\bix}{{\bf x}}
\newcommand{\biy}{{\bf y}}
\renewcommand{\(}{\left(}
\renewcommand{\)}{\right)}

\begin{document}
\title{{\bf Consistent particle systems and duality}}
\author{{
Gioia Carinci$^{\textup{{\tiny(a)}}}$, Cristian Giardin\`a$^{\textup{{\tiny(b)}}}$ and Frank Redig$^{\textup{{\tiny(a)}}}$}
\\
\small $^{\textup{(a)}}$
\small{Delft Institute of Applied Mathematics}\\
\small{Delft University of Technology}\\
{\small Mekelweg 4, 2628 CD Delft, The Netherlands}
\\
\small $^{\textup{(b)}}$
\small{University of Modena and Reggio Emilia}\\
{\small Via G. Campi 213/b, 41125 Modena, Italy}
}
\maketitle
\begin{abstract}
We consider consistent particle systems, which include independent random walkers,
the symmetric exclusion and inclusion processes, as well as  the dual of the KMP model.
Consistent systems are such that the distribution obtained by first evolving $n$ particles and
then removing a particle at random is the same as the one given by a random
removal of a particle at the initial time followed by evolution of the remaining
$n-1$ particles.

In this paper we discuss two main results.
Firstly, we  show that, for reversible systems, the property of consistency is equivalent to
self-duality, thus obtaining a novel probabilistic interpretation of the self-duality
property. Secondly, we  show that consistent particle systems satisfy a set of recursive equations.
This recursions implies that factorial moments of a system with $n$ particles are 
linked to those of a system with $n-1$ particles, thus providing 
substantial information to study the dynamics.
In particular, for a consistent system with absorption, the particle
absorption probabilities satisfy universal recurrence relations.

Since particle systems with absorption are often dual to boundary-driven
non-equilibrium systems, the consistency property implies recurrence relations
for expectations of correlations in non-equilibrium steady states.
We illustrate these relations with several examples.
\end{abstract}
\section{Introduction}
Steady states of non-equilibrium systems, such as systems coupled to reservoirs 
or bulk-driven systems, are of great interest in non-equilibrium statistical mechanics.
However, detailed information such as explicit closed form expression for correlation functions 
are rarely available. In stochastic systems, there are a few integrable systems, such as the symmetric 
and asymmetric exclusion process in dimension one, where one can obtain such 
closed form formulas for all correlation functions in the non-equilibrium steady states 
via matrix product ansatz solution \cite{derrida2}.  Another powerful tool in the analysis 
of non-equilibrium systems is duality, a technique which allows to connect correlation 
functions of order $n$ in the non-equilibrium steady state to a system of $n$ dual
particles. Duality is strictly weaker than integrability, i.e., there are many systems 
which have useful dualities but are not integrable. Duality also allows to study 
time-dependent expectations and time-dependent correlation functions. 
For a correlation function of order $n$, one needs the law of $n$ dual particles.

In particle systems with a duality property \cite{kmp, cggr, demasi2006mathematical,spohn2012large},
whenever we couple the system to appropriately chosen reservoirs,
they are dual to particle systems where the driving reservoirs are replaced by absorbing boundaries.
As a consequence, the computation of correlation functions of the non-equilibrium steady state can
be reduced to the computation of absorption probabilities of dual particles.
These absorption probabilities can usually be obtained explicitly (in closed form) for a single particle, and
in some exceptional cases, such as the symmetric simple exclusion process on a chain \cite{derrida2},
also for many particles. In this paper, we focus on the property of consistency of particle systems
and show that this property provides substantial information to compute the absorption
probabilities.

Consistency, is a property of particle systems generalizing duality and which has a simple probabilistic interpretation.
Intuitively, consistent particle systems are those where the operation of randomly removing a particle commutes with the time-evolution.
This property implies that there are intertwining relations between the dynamics with $n$ and $n-1$ particles, and as a consequence 
simplifying recursion relations.

More in detail, consistent particle systems are defined as permutation invariant particle systems  that additionally  satisfy the following property:
if we marginalize on the first $n-1$ coordinates of the process $\{X^{(n)}(t), \,t\geq 0\}$ describing the
positions of $n$ particles, we exactly obtain (in distribution) the process $\{X^{(n-1)} (t): \geq 0\}$
describing the positions of $n-1$ particles. We require this property to hold for
every number of particles, thus we are actually demanding consistency for a family
of processes, indexed by the number of particles $n$. This  property holds trivially for
independent particles, but  here we are particularly interested in proving and using this property
for systems of interacting particles. The consistency property has already been shown to be
of great relevance in the context of the Kipnis-Marchioro-Presutti (KMP) model
(see \cite{kmp}), where local equilibrium is proved via consistency.

Besides the motivations coming from non-equilibrium statistical physics,
there is clearly a natural probabilistic interest in gaining a better understanding
of the consistency property for interacting particle systems and its relation with (self-) duality. 
Indeed, in this paper we show that consistency and product reversible measures implies 
self-duality, which provides a more probabilistic understanding of self-duality.
As a first step in this direction, we show here that
consistency is indeed equivalent to having a symmetry at the level of the process generator.
More precisely we show that a particle system is consistent if and only if the process generator commutes
with the so-called annihilation operator. This commutation relation formalizes the
property that for consistent particle systems the two operations of ``removing a particle at random''
and ``dynamical evolution'' yield the same result (in distribution) in whatever order they are
performed. The characterization of consistency as a symmetry property expressed by the
vanishing of a commutator allows us to establish a direct link between consistency and
duality \cite{lig}.

A second line of research of this paper is concerned with exploiting the
consistency property to explicitly solve the dynamics, especially in the context of non-equilibrium systems. We will show that
consistency gives rise to a set of recursive equations.
These recursions imply that factorial moments of a system with $n$ particles are
linked to those of a system with $n-1$ particles.
Although the set of recursive equations is not enough to fully
solve the dynamics, it yields substantial simplifications.
In particular, we show the information that is contained
in the recursive equations for a consistent system with absorption
(it is a general property that adding absorbing sites conserves the consistency property).
This leads to a set of recursion relations for the absorption probabilities, which can be iterated to obtain explicit formulas for correlation functions in terms of .
We apply these recursion relations to show universal properties of the non-equilibrium
steady states of boundary-driven systems, such as the exclusion and inclusion process
coupled to boundary reservoirs.

The paper is organized as follows. In section 2 we give the basic definitions
yielding the distinction between the coordinate process, which specifies the positions
of all the particles, and the configuration process, which instead provides the number
of particles at each site. We also recall the definition of the annihilation operator
and its interpretation as a particle removal operator. In section 3 we introduce
the notion of consistent permutation-invariant particle systems and show that this
property is equivalent to a commutation property between the generator of the
configuration process and the annihilation operator. We also write the set
of recursion relation implied by consistency in the general setting.
In section 4, we show the relation between consistency and (self-)duality.
In section 5 we show that consistency is preserved by adding sites where
particles can be absorbed (independently for different particles) and
specify the recursion relations in terms of the absorption probabilities.
Sections 6 and 7 are dedicated to a special class of consistent particle systems,
including the case of an integrable systems (the open symmetric exclusion process
that is solved by matrix product ansatz)
and a non-integrable case (the symmetric inclusion process).
We show in both cases the form taken by consistency equations.

\section{Preliminary definitions}

\subsection{Coordinate process and configuration process}
We consider a system of $n$ particles moving on a countable set of vertices $V$,
with cardinality $|V|$. Their {\em positions}  are denoted by $(x_1,\ldots, x_n)\in V^n$.
The set of functions $g: V^n\to\R$ is denoted by $\caC_n$.
A {\em configuration} of $n$ particles is the $n$-tuple of their positions modulo permutations of labels.
More precisely, for $\bix=(x_1,\ldots, x_n)\in V^n$ the associated configuration is denoted by
\be\label{phi}
\phi(\bix):= \sum_{i=1}^n\delta_{x_i}
\ee
where $\delta_{z}$ is the configuration having only one particle located at $z\in V$, i.e.  for $y\in V$,
$$
(\delta_z)_y = \left\{
\begin{array}{rl}
1 & \text{if } y = z,\\
0 & \text{otherwise}.
\end{array} \right.
$$
We view $\phi$ as a map from $n$-tuples with arbitrary $n$ to configurations, i.e., $\phi: \cup_{n=1}^\infty V^n \to \N^V$.
We then define $\Omega_n$ as the set of configurations
of $n$ particles, i.e.
$$
\Omega_n=\Big\{ \eta\in \Lambda^V:  \|\eta\|=n\Big\}, \qquad \text{with} \qquad \|\eta\|:=  \sum_{x\in V} \eta_x
$$
where $\Lambda \subseteq \N$ is the single-site state space, namely the set of possible occupation numbers for each site. We denote by ${\cal E}_n$ the set of functions $f:\Omega_n\to \R$ and define $\Omega= \{ \eta\in \Lambda^V: \|\eta\|<\infty\}$ the set of finite particle configurations, namely  $\Omega=\cup_{n\in\N}\Omega_n$.

\vskip.2cm
\noindent
In this paper, we shall consider both examples with a finite state space,
such as the partial exclusion processes \cite{lig,schutz1994non} (with a restriction on
the number of particles per site,  i.e. $\Lambda=\{0,\ldots, \alpha\}$
where $\alpha\in \N$ denotes the maximal number of particles per site),
as well as examples with $\Lambda=\mathbb{N}$, such as the
inclusion process \cite{grv} or the independent random walk process 
(see Section \ref{Sect:Ex} for a treatment of these processes).

%

\vskip.2cm
\noindent
In some cases  not all only elements of $V^n$ give rise to allowed configurations of $\Omega$.  For instance, for the partial exclusion process, it is not guaranteed that  $\bix \in V^n$ does not contain more than $\alpha$ particles at any site, i.e.
$\phi(\bix)\in \{0, \ldots, \alpha\}^V$.
For this reason we define the set $V_n$  of $n$-tuples $\bix\in V^n$ such that the associated configuration $\phi(\bix)$  is an element of $\Omega_n$.
For the examples of independent random walkers and inclusion process $V_n=V^n$, whereas, in the case of the partial exclusion process, $V_n\not=V^n$.

\vskip.3cm
\noindent
With these preliminaries we next specify the distinction between
a coordinate process and a configuration process.

\bd[Coordinate process]
We shall call a {\em coordinate process with $n$ particles}, denoted by
 $\{X^{(n)}(t), \,t\geq 0\}$, the stochastic process taking values
in $V_n$ that describes the positions of particles in the course of time. 
Namely, for $i=1,\ldots, n$, the random variable $X_i^{(n)}(t)$  denotes the position of the
$i^{th}$ particle at time $t\geq 0$.   We denote by $\{X(t), t \geq 0\}$ a family of coordinate-processes $\(\{X^{(n)}(t), \,t\geq 0\}, n\in \N\)$, labeled by the number of particles $n\in \N$.
\ed

Throughout this paper we shall restrict to coordinate processes
that are Markov processes.

\bd[Configuration process]
We shall call a {\em configuration process},
denoted by $\{\eta(t), \,t\geq 0\}$,  the stochastic process, 
taking values in $\Omega$ that describes the sites occupancy numbers in the course of time.
Namely, 
for $i\in V$, the random variable $\eta_i(t)$  denotes the  number of particles
at site $i$ at time $t\geq 0$.
\ed

Throughout this paper we shall restrict to configuration processes
that conserve the number of particles, i.e. if the process $\{\eta(t), \,t\geq 0\}$
is started from $\eta\in\Omega_n$ then  $\eta(t)\in\Omega_n$
for all later times $t>0$.
\medskip

A configuration process is naturally induced by a coordinate
process using the map $\phi$ defined in \eqref{phi}. There can be several
coordinate processes whose image under the map  $\phi$
yields the same configuration process. This leads us to the
following definition.

\bd[Compatibility]
A family of coordinate processes $\{ X(t), \,t\geq 0\}$ and a configuration
process $\{\eta(t), \,t\geq 0\}$ are compatible if for all $n\in\N$ the following holds: whenever
$\col{\phi(X_1(0),\ldots, X_n(0))= \eta(0)}$ 
then 
$$
\{\phi(X^{(n)}(t)), \,t\geq 0\}=\{\eta(t), \, t\geq 0 \} 
$$
where the  equality is in distribution.
\ed

Of course it is not guaranteed that, starting from a Markov coordinate
process, the mapping  $\phi$ defined in \eqref{phi} induces
a compatible configuration process that is also Markov.
To further discuss this point we need to introduce the notion of
permutation invariance.
We denote by $\Sigma_n$ the set of permutations of $n$ elements.
Moreover we define the operator  $U_\phi $  mapping functions $f:\Omega_n\to\R$ to functions $U_\phi f\in \caC_n$ via $U_\phi f= f\circ \phi$.

\bd[Permutation invariance]

\noindent
\bi
\item[a)]
A family of coordinate Markov processes $\{ X(t), \,t\geq 0\}$
is said to be {\em permutation-invariant} if,
for every $n\in\N$ and 
permutation $\si\in \Sigma_n$, the
 processes
$\{(X^{(n)}_1(t), \ldots, X^{(n)}_n(t)),\, t\geq 0\}$
and
$\{ (X^{(n)}_{\si(1)}(t), \ldots, X^{(n)}_{\si(n)}(t)),\, t\geq 0\}$
are equal in distribution.
\item[b)] A function $g\in\caC_n$ is said to be
{\em permutation-invariant}
if
$$
g(x_1,\ldots,x_n) = g(x_{\sigma(1)},\ldots,x_{\sigma(n)}) \qquad \text{for all}\quad \si\in \Sigma_n.
$$
Equivalently, a function $g\in\caC_n$ is
{permutation-invariant} if there
exists a function $f:\Omega_n\to\R$ such that
$g= U_\phi f$.

\item[c)]
A probability measure $\mu_n$ on $V_n$ is called {\em permutation-invariant} if, for all $n\in \N$ and  $\si\in \Sigma_n$,
under $\mu_n$, the random vectors $(X_1,\ldots, X_n)$ and $(X_{\si(1)}, \ldots, X_{\si(n)})$ have the same distribution.
\ei
\ed

\noindent
We denote by  $L_n$ the infinitesimal generator of the $n$-particle coordinate  process $\{ X^{(n)}(t), \,t\geq 0\}$, and by  $S_n(t)$ the related  semigroup,  i.e., for $g\in \caC_n$
\[
S_n(t) g(\bix):=\E_{\bix}[g(X^{(n)}(t))]
\]
where $\mathbb{E}_{\bix}$ denotes expectation when the coordinate process
is started from $\bix\in V_n$.

The following lemma shows that to permutation-invariant coordinate processes
one can na\-tu\-rally associate a compatible configuration process enjoying the
Markov property.
\bl\label{conslem}
Let $\{ X(t), \,t\geq 0\}$ be a family of  permutation-invariant coordinate Markov processes with
generators $L_n$, $n\in \N$. Define the operator $\loc$ acting on functions $f:\Omega\to\R$ as
\be\label{calL}
\loc f(\eta):= L_n (U_\phi f)(\bix)
\qquad \text{for all} \quad \eta\in \Omega_n\: \:\text{and}\:\:\bix\in V_n \: : \: \phi(\bix)=\eta
\ee
or, equivalently,
\be\label{calL1}
L_n U_\phi =U_\phi \loc \qquad \text{on} \quad {\cal E}_n.
\ee
Then $\loc$ is the infinitesimal generator of a Markov process $\{\eta(t), t\geq 0\}$ that is a configuration process compatible with  $\{ X(t), \,t\geq 0\}$.

\el
\bpr
For a permutation $\si\in \Sigma_n$ and a function $g\in\caC_n$ we define the operator
\be\label{tsi}
T_\si g(x_1,\ldots, x_n)= g(x_{\si(1)}, \ldots, x_{\si(n)}).
\ee
From the
 permutation invariance
of the family of coordinate Markov processes $\{ X(t), \,t\geq 0\}$ it follows that
\be\label{percom}
[L_n, T_\si]=0 \qquad \text{for all} \quad n\in\N \quad \text{and} \quad \si\in \Sigma_n
\ee
\col{where $[\cdot,\cdot]$ denotes the commutator.}
Let $f:\Omega\to \R$ then, by definition, $U_\phi f$ is a permutation-invariant function.
Hence,  from \eqref{percom} it follows
\be\label{Tsigma}
T_\sigma L_n U_\phi f= L_n T_\sigma U_\phi f= L_n U_\phi f \qquad \text{for all} \quad f\in {\cal E}_n.
\ee
This means that $ L_n U_\phi f$ is permutation-invariant, hence there exists a function $\tilde f: \Omega\to \R$ such that
$$
\col{L_n (U_\phi f)(\bix)=\tilde f(\phi(\bix))=U_\phi \tilde f(\bix)}
$$
namely $L_n U_\phi f= U_\phi \tilde f$. Then it is possible to define the operator $\loc$ acting on functions $f:\Omega \to \R$ such that $\loc  f=\tilde f$, and then \eqref{calL1} is satisfied.
From \eqref{Tsigma} we have that
\be\label{Tsigma2}
T_\sigma S_n(t) U_\phi f= S_n(t) U_\phi f \qquad \text{for all} \quad f\in {\cal E}_n
\ee
and then also $S_n(t) U_\phi f$ is a permutation-invariant function at all times.
Now, if we denote by ${\cal S}(t)$  the semigroup associated to $\loc$, it follows that
\be
{\cal S}(t)f(\eta)=S_n(t)U_\phi f (\bix) \qquad \text{for all} \quad \eta\in \Omega_n\: \:\text{and}\:\:\bix\in V_n \: : \: \phi(\bix)=\eta
\ee
namely
\beq
U_\phi {\cal S}(t)=S_n(t)U_\phi \qquad \text{on}\quad {\cal E}_n
\eeq
From this it follows that $\loc$ is the generator of a Markov process that is a configuration process compatible with  $\{ X(t), \,t\geq 0\}$.
%
%
\epr

%

\subsection{Annihilation operator}

We continue by introducing the operators that remove particles either in the coordinate process
(see \eqref{pil} below) or in the configuration process (see \eqref{anni}).
\bd[Particle removal operators]
For $n\in \N$, $1\leq i \leq n$ we denote by $\pi^{(n)}_i: \caC_{n-1}\to\caC_{n}$ the removal operator of the $i^{th}$ labeled particle, acting on functions
$g\in \caC_{n-1}$ as follows:
\be\label{pil}
(\pi^{(n)}_i g)(x_1,\ldots, x_n)= g(x_1,\ldots, x_{i-1}, x_{i+1}, \ldots, x_n) \qquad \text{for all} \quad x_i\in V
\ee
and we  denote by $\Pi^{(n)}:\caC_{n-1}\to\caC_{n}$ the operator acting  on $g\in \caC_{n-1}$ via
\be\label{pi}
\Pi^{(n)} g= \sum_{i=1}^n \pi^{(n)}_i g.
\ee
We define the ``single-site annihilation operator'' $a$ acting on functions $f:\N_0\to\R$ as
\be\label{a}
a f(n)= \left\{
\begin{array}{ll}
n f(n-1) & \text{if} \:\:n \ge1 \\
 0 & \text{if} \:\: n=0
\end{array}
\right.
\ee
and  the \col{``annihilation operator''}  working on functions $f: \Omega\to\R$ as
\be\label{anni}
\caA f(\eta)=\sum_{x\in V} a_x f(\eta) \quad \text{with} \quad
a_x f(\eta)= \left\{
\begin{array}{ll}
 \eta_x f(\eta-\delta_x) & \text{if} \:\: \eta_x\ge1 \\
 0 & \text{if} \:\: \eta_x=0
\end{array}
\right.
\ee
i.e.  $a_x$ denotes the operator $a$ working on the   variable $\eta_x$.
\ed
We remind the reader here that we are restricting to configurations $\eta\in \Omega$, i.e. with a finite number of particles and
therefore the sum in \eqref{anni} is a finite sum.

The annihilation operator is crucially important in the explanation
of self-dualities for several particle systems \cite{GKR,gkrv}. In particular,
as it will be shown later, the fact that the generator $\caL$ of the
configuration process and the annihilation operator $\caA$
commute, i.e. $[\loc, \caA]= 0$, is enough to obtain a self-duality when the process
has a reversible measure.
 Such a commutation relation for  the configuration process is equivalent to an intertwining relation
between the coordinate process with $n$ particles and the
coordinate process with $n-1$ particles. This equivalence is the object of the next theorem.
\bt\label{intert}
Let $\{ X(t), \,t\geq 0\}$ be a family of coordinate Markov processes  with generators $L_n$, $n\in \N$ and
$\{\eta(t), \,t\geq 0\}$ a compatible Markov configuration process  with generator $\caL$
defined in \eqref{calL}.
Then the following statements are equivalent:
\bi
\item[a)] The generators of the coordinate process with $n$ and $n-1$ particles restricted to permutation
invariant functions are intertwined via $\Pi^{(n)}$, i.e., for every $n\in \N$, and for all $g\in \caC_{n-1}$
permutation-invariant
\be\label{interco}
\col{(L_n\Pi^{(n)})(g)=  (\Pi^{(n)}L_{n-1})(g)}
\ee
\item[b)]  The generator of the configuration process commutes with the total annihilation operator, i.e.
\be\label{anicom}
[\loc, \caA]= 0
\ee
\ei
\et
\bpr
We  first show that
\beq\label{PIn}
\Pi^{(n)} U_\phi = U_\phi \caA \qquad \text{on} \quad {\cal E}_n.
\eeq
This means that
\be\label{crux}
(\Pi^{(n)} (f\circ \phi))(\bix) = \caA f(\eta) \qquad \text{for all} \quad \eta\in \Omega_n\: \:\text{and}\:\:\bix\in V_n \: : \: \phi(\bix)=\eta
\ee
Let $\bix= (x_1, \ldots, x_n)$ and $\eta:=\phi(\bix)=\sum_{i=1}^n \delta_{x_i}$, then  we have
\[
\phi(x_1, \ldots, x_{l-1}, x_{l+1},\ldots, x_n)= \left(\sum_{i=1}^n \delta_{x_i}\right)- \delta_{x_l} = \eta-\delta_{x_l}
\]
As a consequence:
\beq
(\Pi^{(n)} (f\circ \phi))(\bix) &=& \sum_{l=1}^n f( \eta-\delta_{x_l})
\nonumber\\
&=&
\sum_{x\in V} \eta_x f(\eta-\delta_x)
\eeq
where the last step follows because every $x\in V$ is counted exactly $\eta_x$ times in the sum $\sum_{l=1}^n f( \eta-\delta_{x_l})$. This proves \eqref{PIn}.
Suppose now that $[\loc, \caA]=0$, then on ${\cal E}_n$ we have
\[
L_{n-1}  \Pi^{(n)} U_\phi= L_{n-1} U_\phi {\cal A} = U_\phi \loc \caA  = U_\phi \caA\loc = \Pi^{(n)}U_\phi \loc=  \Pi^{(n)} L_n U_\phi
\]
where the equalities follow from \eqref{PIn}, \eqref{calL1} and the commutation relation. Then \eqref{interco} follows since, for all $g\in {\cal C}_n$ permutation-invariant, $g=U_\phi f$ for some $f\in {\cal E}_n$.
The reverse implication is proved analogously.
\epr
\br
The probabilistic interpretation of \eqref{interco} is as follows: if we remove
a randomly chosen particle,  evolve the process, evaluate a permutation-invariant
function at time $t>0$ and finally take expectation, then we can as well
first evolve the process, remove a randomly chosen particle at time $t>0$, evaluate the same permutation-invariant
function and take expectation. In other words,  the operations ``removing a randomly chosen particle''
and ``time evolution in the process followed by expectation'' commute as long as we restrict to permutation-invariant functions.
See also Proposition \ref{CONS}.
\er

\section{Consistency}

In this section we first define the consistency property and 
\col{give some characterization of consistent particle systems.}
We then discuss the implication of consistency in the
form a set of recursive equations.

\subsection{Definition of consistency  and relation to annihilation operator}

\col{Based on Theorem \ref{intert} we define consistent particle systems as follows.}
\bd[Consistency]

\noindent
\bi
\item[a)]
A  family of coordinate Markov processes $\{ X(t), \,t\geq 0\}$
is said to be {\em consistent}  if
 the processes
$\{(X^{(n)}_1(t), \ldots, X^{(n)}_{n-1}(t)), \, t\geq 0\}$
and $\{ (X^{(n-1)}_{1}(t), \ldots, X^{(n-1)}_{n-1}(t)), \, t\geq 0\}$ are equal in distribution for every $n\in\N$.
\item[b)]  A configuration process $\{\eta(t), \; t\ge 0\}$ is said to be {\em consistent} if its generator $\loc$ commutes with the annihilation operator, i.e. $[\loc,\caA]=0$.
\item[c)]
A family of probability measures $\mu_n$ on $V_n$, indexed by $n\in \N$, is called {\em consistent} if for all $n\geq 2$, under $\mu_n$, the
distribution of $(x_1,\ldots, x_{n-1})$ equals $\mu_{n-1}$.
\ei
\ed

\noindent
Let $\{\mu_n, n\in \N\}$ be a consistent collection of probability measures on $V_n, n\in \N$.
If additionally  $\{\mu_n, n\in \N\}$ is permutation-invariant,
then {\em every $m$-dimensional marginal} of $\mu_n$, $m\le n$, coincides with $\mu_{m}$.
A simple example of such a consistent permutation-invariant family is
\[
\mu_n =\frac1{n!} \sum_{\si\in \Sigma_n} \delta_{(x_{\si(1)}, \ldots ,x_{\si(n)})} \qquad \text{for some} \:\: {\bf x}\in V_n
\]
Analogously, if $\{ X(t), \,t\geq 0\}$  is  a family of permutation-invariant coordinate processes,
consistency implies that, for all $n\in \N$,  any $m$-dimensional marginal of $\{ X^{(n)}(t), \,t\geq 0\}$, $m\le n$, is equal in distribution to $\{ X^{(m)}(t), \,t\geq 0\}$, i.e., for all $1\le i_1<\ldots<i_m\le n$,
\be
\{X_{i_1}^{(n)}(t),\ldots,  X_{i_m}^{(n)}(t), \, t\ge 0\}= \{X_1^{(m)}(t),\ldots,  X_m^{(m)}(t), \, t\ge 0\} \qquad \text{in distribution}
\ee
\vskip.3cm
In the following theorem we show that, if a configuration process is consistent then a compatible   coordinate process is also consistent provided that its initial distribution at time $t=0$ is also consistent and permutation invariant.
\bt\label{consthm}
Let $\{ X(t), \,t\geq 0\}$ be a family of  coordinate Markov processes  with generators $L_n$, $n\in \N$,  and let $\{\eta(t), \,t\geq 0\}$ be
a compatible configuration process.
Assume that:
\begin{itemize}
\item[i)]
$\{\eta(t),\; t\ge 0\}$ is consistent;
\item[ii)]
the probability measures $\{\mu_n, n\in \N\}$ on $V_n$, $n\in\N$, form a consistent family which is also permutation-invariant.
\end{itemize}
Then we have consistency of the family of coordinate processes starting from $\{\mu_n, n\in\N\}$, i.e., for all $n\in \N$, $g\in \caC_{n-1}$
permutation-invariant,
\be\label{consi}
\E^{(n)}_{\mu_n} \Big[g(X^{(n)}_1(t), \ldots, X^{(n)}_{n-1}(t))\Big]= \E^{(n-1)}_{\mu_{n-1}} \Big[g(X^{(n)}_1(t), \ldots, X^{(n)}_{n-1}(t))\Big]
\ee
where $\E^n_{\mu_n}$ denotes expectation w.r.t.\ the Markov
process $\{ X^{(n)}(t), \,t\geq 0\}$, started initially with distribution $\mu_n$.
\et
\bpr
Let $g\in \caC_{n-1}$ be a permutation-invariant function, then we have
$\pi^{(n)}_l g= \pi^{(n)}_k g$ for $k,l\in \{1,\ldots, n\}$.
Then, by consistency and permutation invariance of $\mu_n$, we have
\[
\int \pi_l g(x_1,\ldots, x_n) \mu_n(dx_1\ldots dx_n)= \int g(x_1,\ldots, x_{n-1}) \mu_{n-1}(dx_1\ldots dx_{n-1})
\]
 for all $n\in \N$ and all $l\in \{1,\ldots, n\}$.
Therefore, using \eqref{interco}, we have
\beq
\E^{(n)}_{\mu_n} \Big[g(X^{(n)}_1(t), \ldots, X^{(n)}_{n-1}(t))\Big]&=& \E^{(n)}_{\mu_n}\Big[ \pi^{(n)}_{n}(g(X^{(n)}_1(t), \ldots, X^{(n)}_{n-1}(t), X^{(n)}_n(t))\Big]
\nonumber\\
&=&
\int S_n(t)(\pi^{(n)}_{n}g)(x_1,\ldots, x_n) \mu_n (dx_1\ldots dx_n)
\nonumber\\
&=&\frac1n\int\left( S_n(t)\left(\sum_{k=1}^n\pi^{(n)}_{k}g\right)\right)(x_1,\ldots, x_n) \mu_n (dx_1\ldots dx_n)
\nonumber\\
&=&
\frac1n\int S_n(t)(\Pi^{(n)} g)(x_1,\ldots, x_n) \mu_n (dx_1\ldots dx_n)
\nonumber\\
&=&
\frac1n\int (\Pi^{(n)} S_{n-1}(t)g)(x_1,\ldots, x_n) \mu_n (dx_1\ldots dx_n)
\nonumber\\
&=&
\int (S_{n-1}(t)g)(x_1,\ldots, x_{n-1}) \mu_{n-1} (dx_1\ldots dx_n)
\nonumber\\
&=& \E^{(n-1)}_{\mu_{n-1}} \Big[g(X^{(n)}_1(t), \ldots, X^{(n)}_{n-1}(t))\Big].
\eeq
This concludes the proof.
\epr

\subsection{Characterization of consistent configuration processes}
In this section we provide a class of configuration processes exhibithing the consistency property.
We restrict ourselves to models where only one particle jumps at a time. 
We further assume the jump rate to depend only on the number of particles 
hosted by the departure  and arrival sites. 
We consider here generators of the form
\be\label{Gen1}
\loc := \sum_{\{i,j\}\in E}\loc_{i,j}
\ee
where the summmation is over the set $E$ of  (non-oriented) edges $\{i,j\}$ 
of the complete graph with vertices $V$ and
\be\label{Gen11}
 \loc_{i,j}f(\eta)= c_{i,j}(\eta_i,\eta_j)[f(\eta^{i,j})-f(\eta)]+  c_{j,i}(\eta_j,\eta_i)[f(\eta^{j,i})-f(\eta)].
\ee
with  $\eta^{i,j}:=\eta-\delta_i +\delta_j$ and $c_{i,j}: \Lambda\times\Lambda\to \R$ being the  hopping  rate for a particle to jump from site $i$ to site $j$.

We then call $\loc_{i,j}$ the single-edge generator corresponding to the (non-oriented) edge $\{i,j\}$. We will characterize the processes such that
the annihiliation operator commutes with each single edge generator. These processes are then automatically consistent. 

\bt\label{Char}
Let $\{\eta(t), \; t \ge 0\}$ be a configuration process with generator \eqref{Gen1}-\eqref{Gen11}. 
We have that  $[\loc_{i,j}, \caA]=0$  for all $\{i,j\}\in E$ if and only if the rates functions are of the form:
\beq\label{consrate}
c_{i,j}(\kappa,m)=\kappa\cdot \Big(\theta(\{i,j\})\cdot m+\alpha(i,j)\Big)
\eeq
for some $\theta:E\to \R$ and $\alpha:V^2\to \R$.  As a consequence the corresponding configuration process with rates \eqref{consrate} is consistent.
\et
\bpr
We have have to characterize the generators $\loc$ such that for all $i,j\in V$ 
\be\label{lija}
[\loc_{i,j},\caA]=0.
\ee
Recalling that $\caA=\sum_{x\in V}a_x$,
using the assumption that the rates only depend on the number of particles in the departure and arrival sites  \eqref{Gen11}, we have that
 \beq
[\loc_{i,j},a_x]=0 \qquad \text{for all} \quad x\neq i,j 
\eeq
Therefore  we have that \eqref{lija} is satified if and only if  for all $i,j$
\beq
[\loc_{i,j},a_i+a_j]=0
 \eeq
We thus impose
 \beq
 \loc_{i,j}(a_i+a_j)f= \loc_{i,j}f(a_i+a_j)\label{equal}
 \eeq
 for all functions $f:\Lambda\times\Lambda\to \R$ of two variables. 
 One one hand we have that
 \beq
 \loc_{i,j}(a_i+a_j)f(\kappa,m)&=& c_{i,j}(\k,m)\left[(\k-1)f(\k-2,m+1)-\k f(\k-1,m)\right]\nn\\
 &+&c_{j,i}(m,\k)\left[(\k+1) f(\k,m-1)-\k f(\k-1,m)\right]\nn\\
 &+& c_{i,j}(\k,m)\left[(m+1)f(\k-1,m)-m f(\k,m-1)\right]\nn\\
  &+&c_{j,i}(m,\k)\left[(m-1) f(\k+1,m-2)-m f(\k,m-1)\right];
 \eeq
 on the other hand
  \beq
(a_i+a_j) \loc_{i,j}f(\kappa,m)&=& \kappa c_{i,j}(\k-1,m)\left[f(\k-2,m+1)- f(\k-1,m)\right]\nn\\
 &+&\k c_{j,i}(m,\k-1)\left[f(\k,m-1)- f(\k-1,m)\right]\nn\\
 &+& m c_{i,j}(\k,m-1)\left[f(\k-1,m)-f(\k,m-1)\right]\nn\\
  &+&m c_{j,i}(m-1,\k)\left[f(\k+1,m-2)- f(\k,m-1)\right].
 \eeq
 Imposing \eqref{equal} for all $f$ is equivalent to imposing the following conditions for all $m,\k\in \Lambda$:
 \beq
 &&c_{i,j}(\k,m)(\k-1)= c_{i,j}(\k-1,m)\k \label{1}\\
  &&c_{j,i}(m,\k)(m-1)= c_{j,i}(m-1,\k)m\label{2}
  \eeq
  and
  \beq
 && (m+1) c_{i,j}(\k,m)-\k c_{j,i}(m,\k)-\k c_{i,j}(\k,m)= m c_{i,j}(\k,m-1)-\k c_{j,i}(m,\k-1)-\k c_{i,j}(\k-1,m)\nn\\
   && (\k+1) c_{j,i}(m,\k)-m c_{i,j}(\k,m)-m c_{j,i}(m,\k)= \k c_{j,i}(m,\k-1)-m c_{i,j}(\k,m-1)-m c_{j,i}(m-1,\k)\nn\\
  \label{3}
 \eeq
 Iterating the conditions \eqref{1}-\eqref{2} produces
 \beq
&& c_{i,j}(\k,m)=\k \cdot b_{i,j}(m), \qquad b_{i,j}(m):=c_{i,j}(1,m)\label{4}\\
 && c_{j,i}(m,\k)=m \cdot b_{j,i}(\k), \qquad b_{j,i}(\k):=c_{j,i}(1,\k)\label{5}
 \eeq
 for some functions $b_{i,j}, b_{j,i}:\Lambda \to \R$. Inserting \eqref{4} and \eqref{5} in the conditions \eqref{3} we get
 \beq\label{6}
 b_{i,j}(m)-b_{i,j}(m-1)=b_{j,i}(\k)-b_{j,i}(\k-1) \qquad \text{for all} \quad \kappa,m\in \Lambda
 \eeq
 Since the conditions \eqref{6} must be satisfied for all values of $\k$ and $m$, the only possibility is that there exists  a constant $\theta=\theta(\{i,j\})$ such that  
  \beq
 b_{i,j}(m)-b_{i,j}(m-1)=\theta(\{i,j\})=b_{j,i}(\k)-b_{j,i}(\k-1) \qquad \text{for all} \quad \kappa,m\in \Lambda.
 \eeq
 Hence, iterating both for $b_{i,j}$ and for $b_{j,i}$ we deduce that there exist two constants, $\alpha_{i,j}$ and $\alpha_{j,i}$ such that
 \beq\label{7}
 b_{i,j}(m)=\theta(\{i,j\}) m+\alpha(i,j), \qquad  b_{j,i}(\k)=\theta(\{i,j\}) \k+\alpha(j,i)
 \eeq
Then, substituting \eqref{7} in \eqref{4}-\eqref{5} we obtain
  \beq
&& c_{i,j}(\k,m)=\k \cdot \Big(\theta(\{i,j\})\cdot  m+\alpha(i,j)\Big)\\
 && c_{j,i}(m,\k)=m \cdot \Big(\theta(\{i,j\})\cdot \k+\alpha(j,i)\Big)
 \eeq
 that concludes the proof.
\epr
\br
The class of models with rate of the form \eqref{consrate} contains some classical models such as symmetric partial exclusion process, independent random walkers and symmetric inclusion process, see section \ref{Sect:Ex} and \ref{NE} below.
Compared to the class of models with factorized self-duality, described in \cite{sau}, the class of consistent processes introduced in Theorem \ref{Char} is 
larger and in particular contains models which do not have product invariant measures (namely when $\alpha(i,j)\not=\alpha (j,i)$).
\er
\subsection{Recursion relations}

In this section we analyze the consequences of  consistency. More precisely  we prove   recursive relations that  give  the transition probabilities of a system with $m$ particles in terms of the transition probabilities for the system with $m-1$ particles. As a consequence we obtain recursive relations for time dependent factorial moments of consistent  configuration processes.
Specifically, in Theorem \ref{buldog2} below, we show that time dependent factorial moments of order $m$ in a system with $n>m$ particles can be expressed
in terms of $m$-particle transition probabilities. These recursion relations will be  particularly useful when we deal with
systems having absorbing sites (cf.\ sections \ref{5.1} and \ref{TWO}). 
In the reseversible setting they are equivalent with self-duality (cf.\ section \ref{dual-sec}).
\vskip.2cm
\noindent
To prove Theorem \ref{buldog2} , we state two preparatory propositions. 
\bp\label{CONS}
Let $\{\eta(t), \, t\ge 0\}$ a consistent  configuration process on a finite lattice $V$, then, for all $f:\Omega\to \R$ we have
\be\label{conS}
 \sum_{i\in V} \E_{\eta} [\eta_i(t)f(\eta(t)- \delta_i)]=\sum_{i\in V} \eta_i \E_{\eta-\delta_i} [f(\eta(t))]
\ee
\ep
\bpr
From the commutation of the generator of the configuration process $\loc$ with the annihilation operator $\caA$, we obtain the commutation of the semigroup of the configuration process $\caS(t)= e^{t\loc}$ with $\caA$.
We then remark that \eqref{conS} can be written as
\[
\left[\caS(t) (\caA f)\right](\eta) =\left[\caA (\caS(t) f)\right](\eta)
\]
this concludes the proof.
\epr

\noindent
We define now the function
\beq\label{F}
F(\xi,\eta):= \prod_{j\in V}\binom{\eta_j}{\xi_j}, \qquad \xi,\eta\in \Omega
\eeq
and  prove a recursive relation for its expectation. 

\bp\label{buldog1}

Let $\{\eta(t), \, t\ge 0\}$ a consistent  configuration process on a finite  lattice $V$  then, for all $\eta,\xi\in \Omega$ such that $1\le |\xi|< |\eta|$,  we have
\beq\label{cons200}
\E_\eta\left[F(\xi, \eta(t))\right] =\frac 1 {\(|\eta|-|\xi|\)}\cdot \sum_{i\in V} \eta_i \E_{\eta-\delta_i} \left[F(\xi,\eta(t))\right], \qquad \forall \: t\ge 0
\eeq

\ep
\bpr
We  fix $\eta,\xi\in \Omega$ such that $|\xi|\in\{1, \ldots, |\eta|-1\}$ and use \eqref{conS} for the function $f=F(\xi,\cdot)$.
We get
\beq\label{cons1}
\sum_{i\in V} \eta_i \E_{\eta-\delta_i} \left[F(\xi,\eta(t))\right] =\sum_{i\in V} \E_{\eta} \left[\eta_{i}(t)
\binom{\eta_{i}(t)-1}{\xi_i}\prod_{\substack{j\in V\\ j\neq i}}\binom{\eta_{j}(t)}{\xi_j}\right]
\eeq
The  right-hand side of \eqref{cons1} is equal to
\beq\label{o}
 \sum_{i\in V} \E_{\eta} \left[(\eta_i(t)-\xi_i)\cdot \prod_{j\in V}\binom{\eta_{j}(t)}{\xi_j}\right]= (|\eta|-|\xi|)\cdot\E_{\eta} \left[ F(\xi,\eta(t))\right]
\eeq
Combining together \eqref{cons1} and \eqref{o} we obtain \eqref{cons200}.
\epr
\vskip.2cm
\noindent
To formulate the recursion result, we still need some additional notation.
For $n\in \N$, $m\leq n$ we define  $C_{m,n}$, resp. $D_{m,n}$,  to be the set of combinations, resp. dispositions, of $k$ elements chosen  in the set $\{1,\ldots,n\}$, more precisely,
\beq\label{C}
C_{m,n}:=\{(i_1,\ldots, i_m): \: i_j\in \{1,\ldots, n\} \;  \forall \,j\in \{1,\ldots, m\}\: \text{s.t.} \, \:i_1<i_2<\ldots<i_m\}
\eeq
and
\beq\label{D}
D_{m,n}:=\{(i_1,\ldots, i_m): \: i_j\in \{1,\ldots, n\} \;  \forall \,j\in \{1,\ldots, m\}\: \text{s.t.}  \,\:i_j\neq i_k \; \forall \, j\neq k\}
\eeq
Then for $\bix\in V_n$, and $0<m\leq n$ we define
\be\label{cmx}
C_m(\bix)=\{ (x_{i_1}, \ldots, x_{i_m}): \: (i_1,\ldots, i_m)\in C_{m,n}\}
\ee
\bt\label{buldog2}
Let $\{\eta(t), \, t\ge 0\}$ a consistent  configuration process on a finite  lattice $V$. Let $\bix=(x_1,\ldots, x_n)\in V_n$, then, for all $\xi\in \Omega_m$ with $m\in\{1,\ldots, n-1\}$,
\beq\label{Bul}
\E_{\varphi (\bix)}\left[ F(\xi,\eta(t))\right]=\sum_{\biy\in C_m(\bix)} \pee_{\varphi(\biy)}\left(\eta(t)=\xi\right)
\eeq
 where $C_m(\bix)$ is defined in \eqref{cmx}.
\et
\bpr
Let $\eta:=\varphi(\bix)\in \Omega_n$ and $\xi$ as in the hypothesis and define $\kappa:= n-|\xi|=n-m$. From Lemma \ref{buldog1} we have that, if $\kappa \in \{1,\ldots, n-1\}$,
\beq
\E_{\eta}\left[F(\xi,\eta(t))\right]&=&\frac1{n-m}\,\sum_{i\in V} \eta_{i} \, \E_{\eta-\delta_i}\left[F(\xi,\eta(t))\right]\nn\\
&=&\frac1{\kappa}\,\sum_{j_1=1}^n  \E_{\eta-\varphi(x_{j_1})}\left[F(\xi,\eta(t))\right]\label{iter}
\eeq
Now, if $n-1=|\eta-\varphi(x_{j_1})|>m=n-\kappa$, namely, if $\kappa\ge 2$ this can be iterated since
\beq
 \E_{\eta-\varphi(x_{j_1})}\left[F(\xi,\eta(t))\right]=\frac 1 {\kappa-1}\sum_{\substack{j_2=1\\ j_2\neq j_1}}^n  \E_{\eta-\varphi(x_{j_1},x_{j_2})}\left[F(\xi,\eta(t))\right]
\eeq
Let $C_{\kappa,n}$ and $D_{\kappa,n}$ be the sets defined in \eqref{C} and \eqref{D}.
We denote by $J=(j_1,\ldots, j_\kappa)$ an element of $C_{\kappa,n}$ or $D_{\kappa,n}$.
Hence, iterating the argument used in \eqref{iter} $\kappa$ times  we get:
\beq
\E_{\eta}\left[F(\xi,\eta(t))\right]
&=&\frac1{\kappa!}\,\sum_{j_1=1}^n \sum_{\substack{j_2=1\\ j_2\neq j_1}}^n  \ldots  \sum_{\substack{j_\kappa=1\\ j_\kappa\notin \{j_1,\ldots, j_{\kappa-1}\}}}^n\E_{\eta-\varphi(x_{j_1},\ldots, x_{j_\kappa})}\left[F(\xi,\eta(t))\right]\nn\\
&=&\frac1{\kappa!}\,\sum_{J\in D_{\kappa,n}}\E_{\eta-\varphi(x_{j_1},\ldots, x_{j_\kappa})}\left[F(\xi,\eta(t))\right]\nn\\
&=&\sum_{J\in C_{\kappa,n}}\E_{\eta-\varphi(x_{j_1},\ldots, x_{j_\kappa})}\left[F(\xi,\eta(t))\right]\nn\\
&=&\sum_{J\in C_{n-\kappa,n}}\E_{\varphi(x_{j_1},\ldots, x_{j_{n-\kappa}})}\left[F(\xi,\eta(t))\right]\nn\\
&=&\sum_{\biy\in C_{n-\kappa}(\bix)}\E_{\varphi(\biy)}\left[F(\xi,\eta(t))\right]\nn\\
&=&\sum_{\biy\in C_{m}(\bix)}\E_{\varphi(\biy)}\left[F(\xi,\eta(t))\right]\nn\\
\eeq
The theorem follows from the fact that now, for $\biy\in C_{m}(\bix)$,  $|\varphi(\biy)|=m=|\xi|$ and from the observation that
\beq
F(\xi,\eta)= \mathbf 1_{\eta=\xi}, \qquad \text{for} \quad |\eta|=|\xi|
\eeq
so that $\E_{\varphi(\biy)}\left[F(\xi,\eta(t))\right]=\pee_{\varphi(\biy)}\left(\eta(t)=\xi\right)$. This concludes the proof.
\epr

\br The message  of Theorem \ref{buldog2} is the following. Suppose we initialize the con\-fi\-gu\-ra\-tion-process $\{\eta(t), \, t\ge 0\}$ from a configuration $\eta$ with $|\eta|=n$ particles. Now fix $1\le m\le n-1$. Then \eqref{Bul} allows to compute all the factorial moments (of the occupation numbers) of order $m$ in terms of the transition probabilities of the process initialized with $m$ particles. In other words, it is possible to gain information about the system with $n$ particles in terms of the dynamics of $m<n$ particles. Unfortunately the information provided by Theorem \ref{buldog2} is not complete, in the sense that \eqref{Bul} does not give the full distribution of the system with $n$ particles. This is due to the fact that  the factorial-moments of order $n$ are still missing.
\er

We close this section with a specialization of Theorem \ref{buldog2}
that will be useful later. 
\bd\label{assorw}
For a configuration process $\{\eta(t), \;t\geq 0\}$ we define
an associated random walk process  $\{X^{\text{rw}}(t), \, t\ge 0\}$ which is such that when $\eta(0)= \delta_u$ then
$\eta(t)= \delta_{X^{\text{rw}}(t)}$, with $X^{\text{rw}}(0)=u$.
\ed

%
\bc
Let $\{\eta(t), \, t\ge 0\}$ a consistent  configuration process on a finite  lattice $V$ and let $\bix=(x_1,\ldots, x_n)\in V_n$, then, for all $u\in V$ we have
\beq\label{Bul2}
\E_{\varphi (\bix)}\left[\eta_v(t)\right]=\sum_{j=1}^n  \mathbf P_{x_j}(X^{\text{rw}}(t)=v)
\eeq
where $\mathbf P_u$,  is the path space measure of the random walker  $\{X^{\text{rw}}(t), \, t\ge 0\}$ on $V$ starting from  $u\in V$ associated to
the configuration process $\{\eta(t), \; t\geq 0\}$ as defined above in Definition \ref{assorw}.

\ec
\bpr
The result immediately  follows by applying Theorem \ref{buldog2} to the case $\xi=\delta_v$.
\epr

%

\section{Consistency and self-duality}
\label{dual-sec}

In this section we show that consistency implies a form of self-duality whenever a process admits a  strictly-positive reversible measure.
We start by recalling the definition of {\em duality}.
\bd
\label{standard}
Let  $\{Y_t\}_{t\ge0}$,  $\{\widehat{Y}_t\}_{t\ge0}$ be two  Markov processes with state spaces  $\Omega$ and $\widehat{\Omega}$ and $D: \Omega\times \widehat{\Omega}\to\R$ a bounded
measurable function. The processes  $\{Y_t\}_{t\ge0}$,  $\{\widehat{Y}_t\}_{t\ge0}$ are said to be dual with respect to $D$  if
\be\label{standarddualityrelation1}
\E_y \big[D(Y_t, \widehat{y})\big]=\widehat{\mathbb{E}}_{\widehat{y}} \big[D(y, \widehat{Y}_t)\big]\;
\ee
for all $y\in\Omega, \widehat{y}\in \widehat{\Omega}$ and $t>0$. In (\ref{standarddualityrelation1})
$\E_y $ is the expectation with respect to the
 law of the process $\{Y_t\}_{t\ge0}$ started at $y$, while $\widehat{\mathbb{E}}_{\widehat{y}} $
denotes  expectation with respect to the law of the process $\{\widehat{Y}_t\}_{t\ge0}$
initialized at $\widehat{y}$.
\noindent
We say that a process is self-dual when the dual process coincides with the original process.
\ed

It is useful to express the duality property in terms of generators of the two processes. If $L$ denotes the generator
of $\{Y_t\}_{t\ge0}$ and $\widehat L$ denotes the generator of $\{\widehat{Y}_t\}_{t\ge0}$, then (assuming that the duality functions are in the domain of the
generators), the above definition is equivalent to $LD(\cdot,\widehat y)(y)= \widehat L D(y, \cdot)(\widehat y)$ where $L$ acts on the first variable and
$\widehat L$ acts on the second variable.
\vskip.1cm
\noindent
In order to prove our theorem on the relation between consistency and self-duality we recall two general results on self-duality from \cite{gkrv}.
\begin{itemize}
\item[a)] {\em {Trivial} duality function from a reversible measure.}

If a Markov process $\{Y_t:t\geq 0\}$ with countable state-space $\Omega$ has a strictly-positive reversible measure $\nu$, then
the function $D:\Omega\times \Omega \to \R$ given by
\be \label{Cheap}
D(y,\widehat y)= \frac{\delta_{y,\widehat y}}{\nu(y)}
\ee
is a self-duality function.
\item[b)] {\em New duality functions via symmetries.}

If $D:\Omega \times \Omega\to \R$ is a self-duality function and $S$ is a symmetry of $\caL$ (namely an operator acting on functions $f:\Omega\to \R$ commuting with $\caL$, $[S,\caL]=0$) then
$SD$ is a self-duality function, where $SD(y,\widetilde y):=SD(y,\cdot)(\widetilde y)$.
\end{itemize}
\bl\label{lem}
Let $\caA$ denote the annihilation operator defined in \eqref{anni}.
For $\xi, \eta\in \Omega$ denote by $h_\xi(\eta)= \delta_{\xi, \eta}$
the Kronecker delta.
Then we have
\be\label{expa}
(e^{\caA} h_\xi)(\eta)= F(\xi,\eta)
\ee
with $F$ as in \eqref{F}.
\el

\bpr
Recalling the definition of $\caA=\sum_{x\in V} a_x$ where $a_x$ denotes the operator $a$ working on the   variable $\eta_x$ (see \eqref{a}-\eqref{anni}),  in order
to prove \eqref{expa} it is sufficient to show that for all $n,k\in \N_0$
\be\label{expas}
e^a {h}_k(n)=  \binom{n}{k}
\ee
with $h_k(n)=\delta_{k, n}$. For a function $f:\N_0\to\R$
we have
\[
e^a f(n)= \sum_{r=0}^n \binom{n}{r} f(n-r).
\]
Inserting $f= h_k$ gives
\eqref{expas}.
\epr
\vskip.1cm
\noindent
We then obtain the following result.
\bt\label{TD}
Let
$\{\eta(t), \,t\geq 0\}$ be a reversible  configuration process with reversible measure  $\nu$ that is strictly-positive. Then the process is consistent if and only if it is self-dual with self-duality function
\be\label{dualf}
D(\xi,\eta)= \frac{1}{\nu(\xi)}\cdot F(\xi,\eta)
\ee
with $F(\cdot,\cdot)$ as in \eqref{F}.
%
\et
\bpr
We first prove that reversibility and consistency implies that \eqref{dualf} is a self-duality function. From reversibility we know that the function
$$ 
D^{rev}(\xi, \eta)=\frac{\delta_{\xi,\eta}}{\nu(\xi)}
$$
is a self-duality function. If we now act  with $e^{\caA}$ on $D^{rev}$ in the $\eta$-variable
we produce
\eqref{dualf} via Lemma \ref{lem}. This produces a new self-duality function because, by assumption, $\caA$ is a symmetry of $\loc$.
\vskip.1cm
\noindent
Now we prove that reversibility and self-duality with duality function \eqref{dualf} implies consistency.
Fix $\eta\in\Omega_n$ and let $\xi=\eta-\delta_x$ for some $x\in V$, then $\xi\in\Omega_{n-1}$.
Then,
\be\label{bobo}
\nu(\xi) \cdot D(\xi, \eta)=F(\xi,\eta)= \eta_x.
\ee
Let $L_n$ be the operators implicitly defined by \eqref{calL1}. Because $D$ is a self-duality function and $\nu$ is a  strictly-positive function on $\Omega$, satisfying detailed
balance w.r.t. $\loc$, the operator $K: \caE_{n-1}\to \caE_n$
\[
Kf(\eta):=\sum_{\xi\in \Omega_{n-1}} \nu(\xi) D(\xi, \eta) f(\xi)=\sum_{\xi\in \Omega_{n-1}} F(\xi, \eta) f(\xi)
\]
intertwines between $L_n$ and $L_{n-1}$, i.e., $ K(L_{n-1}f)= L_{n} Kf $ for $f\in \caE_{n-1}$ (see \cite{groen} for the connection between duality functions and corresponding intertwining kernel operators).
Now, via \eqref{bobo} we obtain
\[
Kf(\eta)=\sum_x \eta_x f(\eta-\delta_x)= \caA f(\eta)
\]
and  conclude that $ \caA(L_{n-1}f)= L_{n} \caA f$ for $f\in \caE_{n-1}$, which is exactly the commutation property $[\loc, \caA]=0$. The process is thus consistent.
\epr
\section{Consistency for systems with absorbing sites}
In  this section we study consistency for  systems of interacting particles with absorbing sites. These systems are important as they emerge as dual 
processes of boundary-driven non-equilibrium systems connected to reservoirs.
In Section \ref{NE} we will analyze the consequences of consistency
in the context of non-equilibrium systems.

Let $\{\eta(t), \,t\geq 0\}$ denote a configuration process with generator $\loc$.
Let $V^*\supset V$ be a countable set of sites containing $V$, and call $V^{\text{abs}}:=V^*\setminus V$.
We will define a configuration process on $\N^{V^*}$ as follows: inside $V$  particles move according
to the generator $\loc$, and additionally
every particle at site $i$ moves with rate ${ r}(i,j)$ to a site $j\in V^{\text{abs}}$, independently from each other.
Particles arriving at sites $j\in V^{\text{abs}}$ are absorbed and, after absorption, do not move anymore.
We thus obtain what we call the process with absorbing sites $V^{\text{abs}}$ and absorption rates $r(i,j)$.
The generator of this process is  given by
\be\label{absgener}
\loc^{abs} f(\eta) = \loc f(\eta) + {\cal H} f(\eta) \qquad  \quad \text{for} \quad f: \N^{V^*}\to\R
\ee
where $\loc$ only works on the variables $\{\eta_i, i\in V\}$ and where ${\cal H}$ denotes the absorption part of the generator, i.e.,
\be\label{absparto}
{\cal H}f(\eta)=\sum_{i\in V, j\in V^{\text{abs}}} r(i,j) \eta_i \left[f(\eta^{i,j}) -f(\eta)\right] \qquad \text{with} \quad \eta^{i,j}:=\eta-\delta_i +\delta_j
\ee
We can rewrite ${\cal H}$ as follows
\be\label{abspart}
{\cal H}= \sum_{i\in V, j\in V^{\text{abs}}} r(i,j) [a_i a^\dagger_j - a_i a^\dagger_i]
\ee
where $a_i$ is the annihilation operator at site $i$ defined in \eqref{anni} and  $a^\dagger_i$ is the creation operator  defined via
\begin{eqnarray}\label{aiop}
a^\dagger_i f(\eta) &=& f(\eta + \delta_i).\nn
\end{eqnarray}
A process with generator \eqref{absgener} is  called an absorbing extension of the generator $\loc$.
\bd
\vskip.2cm
\noindent
Let $\{X(t), \,t\geq 0\}$ be a   family of coordinate Markov processes on the lattice $V$. Then we define
  its absorbing extension $\{X^{abs}(t), \,t\geq 0\}$ to the lattice $V^*$ as the family of coordinate Markov processes  $\{X^{abs,(n)}(t), \,t\geq 0\}$, $n\in \N$,  on $(V^*)^n$
defined by adding to the jumps of $\{X^{(n)}(t), \,t\geq 0\}$, $n\in \N$, the additional jumps from $i\in V$ to $j\in  V^{\text{abs}}$ at rate $r(i,j)$, that particles perform independently from each other.

\ed
\noindent
Define the set $\Omega^{\text{abs}}$ of configurations on absorbing sites:
\beq
\Omega^{\text{abs}}:=\{ \zeta\in \N_0^{V^{\text{abs}}}: \|\zeta\|<\infty\}.
\eeq
We have the following lemma.
\bl \label{comp} Let $\{ X(t), \,t\geq 0\}$ be a family of  coordinate Markov processes on the lattice $V$  compatible with the configuration process on $\Omega$ having generator $\caL$.  Then its absorbing extension  $\{ X^{abs}(t), \,t\geq 0\}$ to the lattice $V^*$ is compatible with the configuration process on $\Omega^*:=\Omega\times \Omega^{\text{abs}}$ with generator $\loc^{abs}$  in \eqref{absgener}.
\el
\bpr
It follows immediately from the definition of absorbing extensions.
\epr
\vskip.2cm
\noindent
We  then have the following results.
\bl\label{comm} Let $\{\eta(t), \; t\ge 0\}$ be a consistent configuration process on a lattice $V$, then every absorbing extension to a lattice $V^*\supset V$ is  a consistent process.
 \el
\bpr
Let $\caA$ denote the annihilation operator \eqref{anni}.
We want to prove that
\be\label{LA}
[\loc^{abs}, \caA^{abs}]=0  \qquad \quad \text{for}\quad \caA^{abs} f(\eta):= \sum_{i\in V^*} a_i
\ee
with $a_i$ as defined in \eqref{anni}.
Since $[\loc, \caA]=0$ by assumption, and as a consequence $[\loc, \caA^{abs}]=0$,  we only have to prove that
\[
[{\cal H}, \caA^{abs}]=0.
\]
Using \eqref{abspart} and the fact that operators working on variables at different sites commute, we have to show that for all $i,j\in V^*$
\[
[ a_i a^\dagger_j - a_i a^\dagger_i, a_i+a_j]=0
\]
This in turn follows from the commutation relations $[a_i, a_j]=0, [a_i, a^\dagger_j]=\delta_{i,j}$.
\epr

\bt\label{absthm}
Let $\{\eta(t), \; t\ge 0\}$ be a consistent configuration process on the lattice  $V$, and let $\{ X(t), \,t\geq 0\}$ be a family of  coordinate Markov processes   compatible with it.
 Then its  absorbing extension $\{X^{abs}(t), \,t\geq 0\}$ to $V^*\supset V$  is consistent if started from
  a  consistent family  $\{\mu_n, n\in \N\}$ of probability measures on $(V^*)^n$, $n\in \N$,   which is also permutation-invariant. Namely
for all $n\in \N$, $g:(V^*)^n\to \R$
permutation-invariant,
\be\label{consi}
\E^{(n)}_{\mu_n} \Big[g(X^{abs,(n)}_1(t), \ldots, X^{abs,(n)}_{n-1}(t))\Big]= \E^{(n-1)}_{\mu_{n-1}} \Big[g(X^{abs,(n)}_1(t), \ldots, X^{abs,(n)}_{n-1}(t))\Big]
\ee
where $\E^n_{\mu_n}$ denotes expectation w.r.t.\ the Markov
process $\{ X^{abs,(n)}(t), \,t\geq 0\}$, started initially with distribution $\mu_n$.
\et
\bpr
It follows from Lemma \ref{comp}, Lemma \ref{comm} and Theorem \ref{consthm}.
\epr

\subsection{Recursion relations for absorption probabilities}
\label{5.1}

In what  follows we denote by $\pee_{\eta}(\eta(\infty)=\zeta)$
 the probability that eventually $\eta(t)$ settles
in the absorbing configuration $\zeta\in \Omega^{abs}$ starting from the initial configuration $\eta\in \Omega^*$, i.e.,
\be
\pee_{\eta}(\eta(\infty)=\zeta):=\lim_{t\to \infty}\pee_{\eta}(\eta(t)=\zeta).
\ee
Similarly,
\be
\E_\eta[f(\eta(\infty))]=\lim_{t\to \infty}\E_\eta[f(\eta(t))].
\ee

Then, as a consequence of Lemma \ref{comm} and Theorem \ref{buldog2} we have the following   recursion relations for the absorption probabilities.
\bt\label{RR}
Let $\{\eta(t), \, t\ge 0\}$ be a consistent configuration process on a finite  lattice $V^*$ with generator $\loc^{abs}=\loc+{\cal H}$.
Let $\bix=(x_1,\ldots, x_n)\in (V^*)^n$, then, for all $\zeta\in \Omega^{\text{abs}}_m$ with $m\in\{1,\ldots, n-1\}$, we have
\beq\label{rrl}
\E_{\varphi (\bix)}\left[ F(\zeta,\eta(\infty))\right]=\sum_{\biy\in C_m(\bix)} \pee_{\varphi(\biy)}\left(\eta(\infty)=\zeta\right)
\eeq
 with $C_m(\bix)$ as in  \eqref{cmx}.
\et

%
\vskip.3cm
\noindent
The  relations \eqref{rrl} express  combinations of absorption probabilities from the initial configuration $\varphi(\bix)$
in terms of combinations of absorption probabilities from an initial configuration $\eta'$  with  less particles.
  Although these equations are not sufficient to determine the absorption probabilities  in closed form, they are still
  considerably simplifying the problem of computing them as they imply severe restrictions.

\vskip.2cm
\noindent
The following Corollary immediately follows by specializing \eqref{rrl} to the case $\zeta=\delta_v$ for some $v\in  V^{\text{abs}}$.
\bc
Let $\{\eta(t), \, t\ge 0\}$ be a consistent configuration process on a finite  lattice $V^*$ with generator $\loc^{abs}=\loc+{\cal H}$. Let $\bix=(x_1,\ldots, x_n)\in (V^*)^n$, then, for all $v\in V^{\text{abs}}$ we have
\beq\label{Bul1}
\E_{\varphi (\bix)}\left[\eta_v(\infty)\right]=\sum_{j=1}^n  \mathbf P_{x_j}(X^{\text{rw}}(\infty)=v)
\eeq
 where $\mathbf P_u$  is the path-space measure  of the random walk  $\{X^{\text{rw}}(t), \, t\ge 0\}$ on $V^*$ starting from  $u\in V^*$ associated to the configuration process $\{\eta(t), \; t\geq 0\}$ as in Definition \ref{assorw}.

\ec

\subsection{Systems with two absorbing states}\label{TWO}
We now consider the situation in which the system contains only two absorbing states, $|V^{\text{abs}}|=2$, say $V=\{1,\ldots,N\}$ and  $V^{\text{abs}}=\{0,N+1\}$. 
In this case, due to the conservation of particle number, it is sufficient to know the absorption probability in one of the two states, say state 0.  The following proposition is an immediate consequence of  Theorem \ref{RR}.

\bp\label{RR2}
Let $\{\eta(t), \, t\ge 0\}$ be a consistent configuration process on the  lattice $V^*=V\cup V^{\text{abs}}$, $V=\{1,\ldots,N\}$, $V^{\text{abs}}=\{0,N+1\}$, with generator $\loc^{abs}=\loc+{\cal H}$.
Let $\bix=(x_1,\ldots, x_n)\in (V^*)^n$, then, for all  $m\in\{1,\ldots, n-1\}$, we have
\beq\label{rrl2}
 \E_{\varphi (\bix)}\left[ \binom{\eta_0(\infty)}{m}\right]=\sum_{\biy\in C_m(\bix)} \pee_{\varphi(\biy)}\left(\eta(\infty)=m\right)
\eeq
 with $C_m(\bix)$ as in  \eqref{cmx}.
\ep
\bpr This follows by applying Theorem \ref{RR} to the case $\zeta= m\delta_0$,  for  $1\le m\le n-1$.
\epr
\vskip.2cm
\noindent
Notice that \eqref{rrl2} can be viewed as a linear system of $n-1$  independent equations in the $n+1$ variables $\{\pee_\eta \(\eta_0(\infty)= m\), \; m=0,\ldots,n\}$.  Complementing these with the normalization condition  $\sum_{m=0}^{n} \pee_\eta \(\eta_0(\infty)= m\) =1$  we obtain $n$ independent equations. This is still not sufficient to get a closed form expression for the absorption probability (which represent $n+1$ unknowns),
since one independent equation is still missing.

\subsubsection*{Generating function method}
Let    $\{\eta(t), \, t\ge 0\}$ be a consistent process on the  lattice $V^*=V\cup V^{\text{abs}}$, $V=\{1,\ldots,N\}$, $V^{\text{abs}}=\{0,N+1\}$, with generator $\loc^{abs}=\loc+{\cal H}$. We define the  function
\be\label{genplace}
 G(\eta, z):= \E_\eta \left[z^{\eta_0(\infty)}\right], \qquad \eta\in \Omega^*, \: z\ge 0
\ee
The function  $G(\eta,\cdot)$ is the probability generating function of the number of absorbed particles at $0$ starting from the configuration $\eta$. Here we define as usual $G(\eta,0)$ by continuous extension, i.e., $$G(\eta,0):=\lim_{z\to 0}G(\eta,z)= \pee_\eta (\eta_0(\infty)=0)$$ which
is equal to the probability that all the particles in $\eta$ are eventually absorbed at $N+1$.
We then have the following recursion relation.
\bp\label{rec}
Let $\{\eta(t), \, t\ge 0\}$ be a consistent configuration process on the  lattice $V^*=V\cup V^{\text{abs}}$, $V=\{1,\ldots,N\}$, $V^{\text{abs}}=\{0,N+1\}$, with generator $\loc^{abs}=\loc+{\cal H}$. For
all   $\eta\in \Omega^*_n$  we have
\be\label{bonan}
(1-z)G'(\eta,z)  + nG(\eta, z) = \sum_{i\in V^*} \eta_i G(\eta-\delta_i, z)
\ee
and, as a consequence,
\be\label{dualrecu}
G(\eta,z)= (1-z)^n G(\eta,0) + (1-z)^n\sum_{i\in V^*}  \eta_i \int_0^z \frac{1}{(1-u)^{n+1}}\, G(\eta-\delta_i, u) du
\ee
\ep
\bpr
Because of commutation of the generator of the absorbing system with the annihilation operator \eqref{anni} we can write
\[
\lim_{t\to\infty}\E_\eta [(\caA z^{\eta_0}) (t)]= \sum_{i\in V^*} \eta_i G(\eta-\delta_i, z)
\]
which leads
to
\[
\E_{\eta}[\eta_0(\infty) z^{\eta_0(\infty) -1}] + \E_{\eta}\left[(n-\eta_0(\infty)) z^{\eta_0(\infty)}\right]= \sum_{i\in V^*} \eta_i G(\eta-\delta_i, z)
\]
which gives \eqref{bonan}.
We can now ``integrate'' the recursion \eqref{bonan} as follows. Putting $G(\eta, z)=(1-z)^n H(\eta, z)$ and substituting in \eqref{bonan}
we find,
\[
(1-z)^{n+1} H'(\eta,z)= \sum_{i\in V^*}  \eta_i G(\eta-\delta_i, z).
\]
Noticing that $H(\eta,0)= G(\eta,0)$, by  integrating we obtain \eqref{dualrecu}.
\epr

\noindent
The recursion \eqref{dualrecu} can be iterated until we are left with one particle configurations for which, in ideal situations, the generating function can be computed, as we will see in the last section. It is crucial then to obtain a formula for the function
\beq
{\gee}(\eta,z):=G(\eta,z)-G^{\text{irw}}(\eta,z)
\eeq
that is the difference between the generating function of our process $\{\eta(t), \, t\ge 0\}$ and the generating function $G^{\text{irw}}(\eta,z)$  of the auxiliary process $\{\eta^{\text{irw}}(t), \, t\ge0\}$ of independent random walkers defined as follows. Let
 $\{X^{\text{rw}}(t), \, t\ge 0\}$ be the random walker on $V^*$ associated to  the configuration process $\{\eta(t), \, t\ge 0\}$ as in Definition \ref{assorw}.  Now let $\{X^{\text{irw}}(t), \, t\ge 0\}$ be the family of coordinate processes $\{X^{\text{irw},(n)}(t), \, t\ge 0\}$, $n\ge 1$ on $(V^*)^n$
 whose coordinates are $n$ independent copies of $X^{\text{rw}}(t)$:
 $$
 X^{\text{irw},(n)}(t)=(X^{\text{rw}}_1(t),\ldots, X^{\text{rw}}_n(t)), \qquad t\ge 0.
 $$
Then we define  $\{\eta^{\text{irw}}(t),  \, t\ge0\}$ as the configuration  process compatible with $\{X^{\text{irw}}(t), \, t\ge 0\}$  and  we denote by $\pee^{\text{irw}}_\eta$ the related path-space measure conditioned to $\eta^{\text{irw}}(0)=\eta$.
\vskip.2cm
\noindent
It is clear  that for $|\eta|=1$, $\gee(\eta,z)=0$. In the next theorem we obtain a formula for the difference function $\gee(\eta,\cdot)$ when $|\eta|\ge 2$.
\bt\label{GEE}
Let $\{\eta(t), \, t\ge 0\}$ be a consistent configuration process on the  lattice $V^*=V\cup V^{\text{abs}}$, $V=\{1,\ldots,N\}$, $V^{\text{abs}}=\{0,N+1\}$, with generator $\loc^{abs}=\loc+{\cal H}$. For all $\bix\in(V^*)^n$ we have
\beq\label{gee}
\gee(\phi(\bix),z)=\sum_{\kappa=2}^n z^{n-\kappa}(1-z)^\kappa \sum_{\biy\in C_\kappa(\bix)} \gee(\phi(\biy),0),
\eeq
with
\beq
\gee(\eta,0)=\pee_{\eta}(\eta_0(\infty)=0)-\pee^{\text{irw}}_{\eta}(\eta_0(\infty)=0)
\eeq
and $C_\kappa(\bix)$ as in  \eqref{cmx}.
\et
\bpr We proceed by induction on $n$. We first consider the case  $n=2$. In this case
\beq
\sum_{i\in V^*}\eta_iG(\eta-\delta_i, z)=\sum_{i\in V^*}\eta_iG^{\text{irw}}(\eta-\delta_i, z)
\eeq
because $G(\zeta, z)$ and $G^{\text{irw}}(\zeta, z)$ coincide on configuration $\zeta$ with one particle.
Therefore,
\beq
\sum_{i\in V^*}\eta_i\gee(\eta-\delta_i, z)=0
\eeq
and we obtain from Proposition \ref{rec},
\beq
\gee(\phi(\bix),z)=(1-z)^2 \gee(\phi(\bix),0)
\eeq
which is \eqref{gee} for $n=2$ . Now we assume that \eqref{gee} is true for $n-1$ and prove the induction step.
First of all we notice that we can rewrite
\beq
\sum_{i\in V^*}(\phi(\bix))_iG(\phi(\bix)-\delta_i,z)=\sum_{j=1}^n G(\phi(\bix)-\delta_{x_j},z).
\eeq
Using \eqref{dualrecu}  and the induction hypothesis we have that
\beq
\frac 1 {(1-z)^n}\cdot \gee(\phi(\bix),z)-  \gee(\phi(\bix),0)
&&= \sum_{j=1}^n \int_0^z \frac{1}{(1-u)^{n+1}} \,  \gee(\phi(\bix)-\delta_{x_j},u)\,du\nn\\
&&= \sum_{j=1}^n  \sum_{\kappa=2}^{n-1} \sum_{\biy\in C_\kappa(\bix\setminus x_j)} \gee(\phi(\biy),0) \int_0^z \frac{u^{n-\kappa-1}}{(1-u)^{n-\kappa+1}}   \, du
\nn
\eeq
Calling $m=n-\kappa-1$, we have
\beq
&&\int_0^z\frac{u^{m}}{(1-u)^{m+2}}   \, du= \int_{1-z}^1 \frac 1 {v^2}\(\frac 1 v -1\)^m\, dv\nn\\
&&=\sum_{k=0}^m \binom{m}{k}(-1)^{m-k} \int_{1-z}^1 \frac {dv} {v^{k+2}}\nn\\
&&=\sum_{k=0}^m \binom{m}{k}(-1)^{m-k} \frac 1 {k+1} \(\frac 1 {(1-z)^{k+1}}-1\)\nn\\
&&=\frac 1 {m+1}\cdot \sum_{k=0}^m \binom{m+1}{k+1}(-1)^{(m+1)-(k+1)}  \(\frac 1 {(1-z)^{k+1}}-1\)\nn\\
&&=\frac 1 {m+1}\cdot \sum_{l=0}^{m+1} \binom{m+1}{l}(-1)^{(m+1)-l}  \(\frac 1 {(1-z)^{l}}-1\)=\frac 1{m+1}\(\frac z {1-z}\)^{m+1}
\eeq
Hence
\beq\label{QUA}
&&\frac 1 {(1-z)^n}\cdot \gee(\phi(\bix),z)-  \gee(\phi(\bix),0) \nn\\
&&=   \sum_{\kappa=2}^{n-1}\(\frac z {1-z}\)^{n-\kappa}\frac 1{n-\kappa}
\sum_{j=1}^n \sum_{\biy\in C_\kappa(\bix\setminus x_j)} \gee(\phi(\biy),0)
\eeq
Consider now
\beq
&&\sum_{j_0=1}^n \sum_{\biy\in C_\kappa(\bix\setminus x_{j_0})} \gee(\phi(\biy),0) =\frac 1 {\kappa!}\sum_{j_0=1}^n \sum_{\substack{j_1=1\\j_1\neq j_0}}^n\ldots \sum_{\substack{j_\kappa=1\\j_\kappa\neq j_0, \ldots j_{\kappa-1}}}^n \gee(\phi(x_{j_1},\ldots,x_{j_\kappa}),0)\nn\\
&&=\frac 1 {\kappa!} \sum_{\substack{j_1=1}}^n\sum_{\substack{j_2=1\\j_2\neq j_1}}^n\ldots \sum_{\substack{j_\kappa=1\\j_\kappa\neq j_1, \ldots j_{\kappa-1}}}^n \sum_{\substack{j_0=1\\j_0\neq j_1, \ldots j_{\kappa}}}^n \gee(\phi(x_{j_1},\ldots,x_{j_\kappa}),0)\nn\\
&&=\frac {n-\kappa} {\kappa!} \sum_{\substack{j_1=1}}^n\sum_{\substack{j_2=1\\j_2\neq j_1}}^n\ldots \sum_{\substack{j_\kappa=1\\j_\kappa\neq j_1, \ldots j_{\kappa-1}}}^n \gee(\phi(x_{j_1},\ldots,x_{j_\kappa}),0)\nn\\
&&= (n-\kappa)\sum_{\biy\in C_\kappa(\bix)} \gee(\phi(\biy),0)\nn
\eeq
Thus, from \eqref{QUA} we get
\beq\label{QUA2}
&&\frac 1 {(1-z)^n}\cdot \gee(\phi(\bix),z)-  \gee(\phi(\bix),0) =   \sum_{\kappa=2}^{n-1}\(\frac z {1-z}\)^{n-\kappa}\sum_{\biy\in C_\kappa(\bix)} \gee(\phi(\biy),0)
\eeq
that concludes the proof.
\epr

\section{A class of consistent processes}\label{Sect:Ex}

In this section  we consider a natural class of consistent configuration processes $\{\eta(t), \, t\ge 0\}$. These can be obtained,  with a certain choice of the parameters, as particular cases of the more general class of processes produced in the characterization Theorem \ref{Char}.  These processes  do not constitute  the entire class of   processes exibithing the consistency property. Nevertheless they are paradigmatic examples within this class, as they are well known in the literature. The processes we consider are of three types: partial exclusion processes, inclusion processes and independent random walkers. For the sake of synthesis we formally define a unique  generator  
 and let this be  parametrized by a constant $\theta\in \R$ tuning the  attractive or repulsive nature of particle-interaction.  The generator is given by
\be\label{Gen}
\loc_\theta f(\eta):= \sum_{i,j\in V}p(i,j)\, \eta_i(1+\theta\eta_{j})[f(\eta^{i,j})-f(\eta)]
\ee
acting on  functions $f:\Omega_{\theta} \rightarrow \R$, with $\Omega_\theta$ to be defined. Here $p:V\times V\to \R$ is a symmetric function: $p(i,j)=p(j,i)$ and $\eta^{i,j}:=\eta-\delta_i +\delta_j$. Moreover, for simplicity, we assume that it is finite range, i.e., $p(i,j)=0$ if $|i-j|>R$ for some $R\geq 1$.

Remark that if the configuration contains only one particle then this particle moves according to a continuous-time random walk
jumping from $i$ to $j$ at rate $p(i.j)$, which is not depending on the interaction parameter $\theta$.
The parameter  $\theta \in \R$ can also be negative, under the condition that $\frac 1 {|\theta|}$ is integer: $\theta\in\R^+\cup \{\alpha<0: \: -1/\alpha\in \N\}$ and,  according to its sign  we recover one of  three cases: the partial symmetric exclusion-process SEP$(1/|\theta|)$, the symmetric inclusion process SIP$(1/|\theta|)$, and the independent-random-walkers process IRW:
\be
\locb_\theta= \left\{
\begin{array}{ll}
\loc^{\text{IRW}} & \text{for } \theta=0 \\
|\theta|\loc^{\text{SIP}(1/|\theta|)}  & \text{for } \theta>0 \\
|\theta|\loc^{\text{SEP}(1/{|\theta|})}  & \text{for } \theta<0, \frac 1 {|\theta|}\in \N \\
\end{array}
\right.
\ee
Also the state space $\Omega_\theta$ where configurations $\eta$ live changes according to the choice of $\theta$, we have:
\be\label{space}
\Omega_\theta= \Lambda_\theta^{V}, \qquad \text{with} \quad \Lambda_\theta= \left\{
\begin{array}{ll}
 \N & \text{for } \theta\ge 0 \\
\{1,2,\ldots, \tfrac 1 {|\theta|}\} & \text{for } \theta<0, \frac 1 {|\theta|}\in \N \\
\end{array}
\right.
\ee
These processes have been introduced in \cite{gkrv} and broadly studied due to their algebraic properties. In \cite{gkrv} it has been proved that
$[\loc,\caA]=0$ where $\caA$ is the annihilation operator defined in \eqref{anni}.
Notice that it is possible to rewrite the generator as
\be\label{Gen0}
\loc_\theta f(\eta):= \sum_{\{i,j\}\in E}p(\{i,j\})\,\loc_{\theta,i,j}f(\eta)
\ee
where now the summmation is over the edges $\{i,j\}\in E$ of the complete graph with vertices $V$ and
\be\label{Gen00}
 \loc_{\theta,i,j}f(\eta)= \eta_i(1+\theta\eta_{j})[f(\eta^{i,j})-f(\eta)]+  \eta_j(1+\theta\eta_{i})[f(\eta^{j,i})-f(\eta)].
\ee
It is possible to verify that the commutation relation with the annihilation operator holds true at the level of each bond, namely $[\loc_{\theta,i,j},\caA]=0$ for all $i,j\in V$.
\vskip.2cm
\noindent
The processes   admit an infinite family of reversible homogeneous product measures $\nu_{\rho,\theta}$ on $\Omega_\theta$ labelled by the density parameter  $\rho:=\langle \eta_i \rangle_{\nu_{\rho,\theta}} >0$, $i\in \Z$,
with marginals
\be\label{statmeas}
\nu_{\rho,\theta}(\eta_i=n)=\frac {\rho^n} {n!} \cdot\left\{
\begin{array}{lll}
e^{-\rho} & \text{for } \theta=0& \text{Pois}(\rho)\\
&\\
(1+\theta \rho)^{-n-\tfrac 1 \theta} \cdot |\theta|^n (1/|\theta|)^{(n)} & \text{for } \theta>0 & \text{NegBin}(\tfrac1\theta, \tfrac{\theta\rho}{1+\theta\rho})\\
&\\
(1+\theta \rho)^{-n-\tfrac 1 \theta} \cdot |\theta|^n (1/|\theta|)_{n} & \text{for } \theta<0& \text{Bin}(\tfrac 1 {|\theta|}, |\theta|\rho)\\
\end{array}
\right.
\ee
for all $i \in\Z$, $n\in \Lambda_\theta$. Hence, from Theorem \ref{TD}, the processes are self-dual with duality functions of the form $D(\xi,\eta)=F(\xi,\eta)/\nu_\rho(\xi)$ with $F$ as in \eqref{F} (modulo a factor that depends on the total number of particles $|\xi|$). More precisely, the self-duality functions are given by
\be\label{D0}
D_\theta(\xi,\eta)= \prod_{i\in V} d_\theta(\xi_i,\eta_i), \quad \text{with} \quad
d_\theta(k,n)= \frac{n!}{(n-k)!}\cdot \left\{
\begin{array}{ll}
1 & \text{for } \theta= 0 \\
& \\
\frac{1}{\(1/{|\theta|}\)^{(k)}} & \text{for } \theta> 0 \\
& \\
\frac{1} {\(1/{|\theta|}\)_{k}}& \text{for } \theta<0, \frac 1 {|\theta|}\in \N \\
\end{array}
\right.
\ee
where $(a)^{(n)}$ and $(a)_n$ are the Pochhammer symbols for  rising and falling factorials:
\beq
(a)^{(n)}:= \frac{\Gamma(a+n)}{\Gamma(a)}\qquad \text{and} \qquad (a)_n:=\frac{\Gamma(a+1)}{\Gamma(a+1-n)}
\eeq
From consistency we know that for  the processes with generator \eqref{Gen} Theorem \ref{buldog2} holds true.
\subsection{Adding absorbing sites}
Let now $V^*=V\cup V^{\text{abs}}\subseteq \Z^d$ with $V^{\text{abs}}$ a set with absorbing sites.
Define $\loc_\theta^{\text{abs}}:=\loc_\theta+{\cal H}$ with $\loc_\theta$ as in \eqref{Gen} and ${\cal H}$ as in \eqref{absparto}:
\be\label{GenAb}
\loc^{\text{abs}}_\theta f(\eta):= \sum_{i,j\in V}p(i,j)\, \eta_i(1+\theta\eta_{j})[f(\eta^{i,j})-f(\eta)]+\sum_{i\in V, j\in V^{\text{abs}}} { r}(i,j)\, \eta_i \left[f(\eta^{i,j}) -f(\eta)\right].
\ee
From Lemma \ref{comm} we know that $[\loc_\theta^{\text{abs}},\caA^{\text{abs}}]=0$, and as a consequence  for these processes we can apply Theorem \ref{RR}
when studying the absorption probabilities.

\subsection{Adding reservoirs}

We  now add to the bulk generator \eqref{Gen} additional terms describing  the action of external reservoirs each acting on one of the sites of the lattice $V$, or, eventually, on a sample of selected sites.
The generator  has the following form:
\begin{eqnarray}\label{GenRes}
{\loc_\theta^{\text{res}}} f(\eta)&=& \sum_{i,j\in V}p(i,j)\, \eta_i(1+\theta\eta_{j})[f(\eta^{i,j})-f(\eta)]\\
&+&\sum_{i\in V^{\text{ext}}}c(i)\left\{\rho_i (1+\theta \eta_i)[f(\eta+\delta_i)-f(\eta)]+(1+\theta \rho_i)  \eta_i [f(\eta-\delta_i)-f(\eta)]\right\}.\nonumber
\end{eqnarray}
for some $V^{\text{ext}}\subseteq V$ and  $\rho,c:V^{\text{ext}}\to [0, +\infty)$.
Particles are injected and removed from the system through the sites $i\in V^{\text{ext}}$ (which we call the sites coupled to reservoirs) at suitable rates. Here $c(i)$ is the global rate at which the external reservoirs acts on the site $i$ and $\rho_i$ is the density imposed by the reservoir on that site. This means that for each reservoir term, the stationary distribution
for the particle number of that site is exactly $\nu_{\rho_i,\theta}$. In particular if all $\rho_i=\rho$ are equal, then the reversible measure is $\nu_{\rho,\theta}$
which is the equilibrium setting. In all other cases, we have that the system settles in a non-equilibrium steady state with non-trivial correlations as soon as
$\theta\not=0$. In order to avoid uninteresting degenerate cases, we will always choose $p(i,j)$ and $c(i)$ in such a way that there exists a unique stationary measure, which we denote by $ \mu_\theta^{\text{st}} $.
The reservoirs rates have been chosen in such a way to preserve a duality property. Indeed, in \cite{cggr} we proved that the process with generator  \eqref{GenRes} is dual to a system with absorbing sites, i.e.,  generated by \eqref{GenAb}. Here the set of absorbing site $V^{\text{abs}}$ is  a ``copy'' of $V^{\text{ext}}$, i.e., such that
$|V^{\text{abs}}|=|V^{\text{ext}}|$ and the rates  $r(\cdot,\cdot)$ of the absorbing part of the generator \eqref{GenAb} are given by:
$$
r(i,j)=\mathbf 1_{i\in V^{\text{ext}}, j=i^*}\cdot c(i)
$$
where $*$ is a bijection $*:V^{\text{ext}}\to V^{\text{abs}}$ that assigns to each site $i \in V^{\text{ext}}$ a corresponding absorbing site $i^*\in V^{\text{abs}}$. From now on we will denote by $\{\eta(t), \, t\ge 0\}$ the process with reservoirs with generator \eqref{GenRes}  and state space $\Omega_\theta=\Lambda_\theta^V$ (as in \eqref{space}) and by $\{\xi(t), \, t\ge 0\}$ the dual process with absorbing sites and with generator:
\be\label{GenDual}
\loc_\theta^{\text{Dual}} f(\xi):= \sum_{i,j\in V}p(i,j)\, \xi_i(1+\theta\xi_{j})[f(\xi^{i,j})-f(\xi)]+\sum_{i\in V^{\text{ext}}} c(i)\, \xi_i \left[f(\xi^{i,i^*}) -f(\xi)\right]
\ee
and state space $\Omega^*_\theta=\Lambda_\theta^{V^*}$, $V^*=V\cup V^{\text{abs}}$. The duality function $\widehat D_\theta:\Omega_\theta^*\times \Omega_\theta\to \R$ between the two processes is the following:
\be
\widehat D_\theta(\xi,\eta)=\prod_{i\in V} d_\theta(\xi_i,\eta_i) \cdot \prod_{i\in V^{\text{abs}}}\rho_i^{\xi_i}= D_\theta(\xi,\eta)\cdot \prod_{i\in V^{\text{abs}}}\rho_i^{\xi_i}
\ee
with $d$ is as in \eqref{D0}  and where we extend the definition of density function  $\rho$ to the absorbing sites, $\rho:V^{\text{ext}}\cup V^{\text{abs}}$, by identifying: $\rho_i=\rho_{i^*}$. See \cite{cggr} for a proof of the duality relation.
\vskip.2cm
\noindent
An important consequence of this duality property is the information that it gives about the non-equilibrium stationary measure $\mu_\theta^{\text{st}}$ of the process generated by \eqref{GenRes}. Indeed, the expectations of duality functions under $\mu_\theta^{\text{st}}$  can be expressed  in terms of the absorption probabilities of the dual process:
\be\label{corform}
\int \widehat D_\theta(\xi, \eta) \mu_\theta^{\text{st}} (d\eta) = \E_\xi\bigg[\prod_{i\in V^{\text{abs}}} \rho_i^{\xi_i(\infty)}\bigg]
\ee
Notice that now, thanks to  the consistency property of the dual process $\{\xi(t), \, t\ge 0\}$, we can apply Theorem \ref{RR} to get   information about the moments in \eqref{corform}. In the next section  we specialize  to the case of two reservoirs (and then two dual absorbing sites) for which we can use the generating-function method developed in  Section \ref{TWO}.
\vskip.2cm
\noindent
In the case $\theta=0$ we have independent walkers, and the non-equilbrium steady state is a product of Poisson measures with
parameter given by the local density.
For general $\theta$, if there is only one dual particle, both the duality function and the random walk of that dual particle is not depending on the interaction parameter $\theta$.
As a consequence we have for all $x\in V$
\beq\label{onedualpart}
\E_{\mu_\theta^{\text{st}}} (\eta_x)=\int \widehat D_\theta(\delta_x, \eta) \mu_\theta^{\text{st}} (d\eta) &=&
\int \widehat D_0(\delta_x, \eta) \mu_0^{\text{st}} (d\eta)  \nonumber\\
&=&
\E_{\mu_0^{\text{st}}} (\eta_x)=\sum_{y\in  V^{\text{abs}}} \rho_y q(x,y)
\eeq
where $q(x,y)$ is the probability for the walker starting at $x$ to be absorbed at $y$ eventually:
\beq\label{ABS}
q(x,y):={\bf P}_x(X^{\text{rw}}(\infty)=y)
\eeq
Finally, if the dual configuration $\xi=\sum_{i=1}^n \delta_{x_i}$ where all $x_i$ are mutually different elements of $V$, then
\be\label{prod}
 \widehat{D}_\theta(\xi, \eta) = \widehat{ D}_0(\xi, \eta) =\prod_{i=1}^n \eta_{x_i}
\ee

\subsection{Instantaneous thermalization models}
Another class of models sharing the consistency property is the class of processes  obtained as the instantaneous thermalization limit of the processes defined in \eqref{Gen} which we now briefly describe. An important example in this class is the dual KMP model, see \cite{cggr,kmp}. An instantaneous thermalization process gives rise, for each couple of nearest neighbouring sites, to an
instantaneous redistribution of the total number of particles.  For each
bond, the total number of particles in that bond $\eta_i+\eta_j$ is  redistributed according
to the stationary measure of the original process at equilibrium on that bond, conditioned to the
conservation of $\eta_i+\eta_j$. The generators of the instantaneous thermalization processes we consider here is given by
\beq\label{therm}
\loc^{\text{th}}_\theta f(\eta):=\sum_{\{i,j\}\in E} p(i,j) \loc^{\text{th}}_{\theta,i,j} f(\eta)\qquad \text{with} \qquad \loc_{\theta,i,j}^{\text{th}}:=\lim_{t\to \infty}\(e^{t\loc_{\theta,i,j}}-\mathbf 1\)f
\eeq
with $\loc_{\theta,i,j}$ as in \eqref{Gen00}. Since we have that the commutation $[\loc_{\theta,i,j},\caA]=0$ for all
 $i,j\in V$, it follows that also
$\loc^{\text{th}}_{\theta,i,j}$ commutes with $\caA$ and then the thermalized models \eqref{therm} are also consistent, namely  $[\loc^{\text{th}}_\theta,\caA]=0$.
A more explicit  expression of the generator is given by:
\beq\label{therm0}
 \loc_{\theta,i,j}^{\text{th}}:=\sum_{m=0}^{\eta_i+\eta_j}[f(\eta^{i,j,m})-f(\eta)]\; \cdot \bar \nu_\theta(m \;| \, \eta_i+\eta_j)
 \eeq
where
\be
\eta_k^{i,j,m}:=
\left\{
\begin{array}{ll}
\eta_k & \text{for} \:k\neq i,j\\
m &\text{for} \: k=i\\
\eta_i+\eta_j-m & \text{for} \: k=j\\
\end{array}
\right.
\ee
and
$\bar \nu_\theta(m \;| \, M):=\nu_{\theta,\rho}(\eta_i=m\;|\: \eta_i+\eta_j=M)$ with $\nu_{\theta,\rho}$ the reversible measure defined in \eqref{statmeas}.
The process \eqref{therm0} is self-dual with duality function \eqref{D0} (see Section 5 of \cite{cggr}).
\vskip.2cm
\noindent

Also for the instantaneous thermalization models it is possible to add absorbing boundaries in such a way to preserve the consistency property, and also in this case there is duality with the system with absorbing boundaries and a system with reservoirs (if the action of reservoirs is properly chosen, see \cite{cggr}).
\vskip.2cm
\noindent
In the next section  we use consistency to obtain an  expression  for the $n$-points stationary correlation function for the system with reservoirs \eqref{GenRes} in a specific setting. This result can be easily extended to the thermalized models with reservoirs, since, as we have seen here, they share the same commutation property and as a consequence consistency property.

\section{Correlation functions in non-equilibrium steady states}\label{NE}

Here we consider processes in the class of models introduced in the previous section, and, as in  Section \ref{TWO}, we restrict  to the case where $V=\{1,\ldots,N\}$, with $V^{\text{ext}}=\{1,N\}$ and $V^{\text{abs}}=\{0,N+1\}$.
In the spirit of the previous section we assign to each site coupled to a reservoir an absorbing site, namely we say that $1^*=0$ and $N^*=N+1$. Then we
denote by $\{\eta(t), t\ge 0\}$ on $\Omega_\theta$ the process with generator:
\begin{eqnarray}\label{GenRes1}
{\loc_\theta^{\text{res}}}f(\eta)&=& \sum_{i,j\in V}p(i,j)\, \eta_i(1+\theta\eta_{j})[f(\eta^{i,j})-f(\eta)]\\
&+&\cll\left\{\rll(1+\theta \eta_1)[f(\eta+\delta_1)-f(\eta)]+(1+\theta \rll)  \eta_1 [f(\eta-\delta_1)-f(\eta)]\right\}.\nonumber\\
&+&\crr\left\{\rrr (1+\theta \eta_{N})[f(\eta+\delta_{N})-f(\eta)]+(1+\theta \rrr)  \eta_{N} [f(\eta-\delta_{N})-f(\eta)]\right\}.\nonumber
\end{eqnarray}
for some $\cll,\crr,\rll,\rrr\ge 0$, and by $\{\xi(t), \, t\ge 0\}$ its dual with state space $\Omega_\theta^*$ and generator:
\beq\label{GenDual1}
\loc_\theta^{\text{Dual}} f(\xi)&:=& \sum_{i,j\in V}p(i,j)\, \xi_i(1+\theta\xi_{j})[f(\xi^{i,j})-f(\xi)]\nn\\
&+& \cll\, \xi_1 \left[f(\xi^{1,0}) -f(\xi)\right]+ \crr\, \xi_N \left[f(\xi^{N,N+1}) -f(\xi)\right]
\eeq
with duality function given by
\be
\widehat D_\theta(\xi,\eta)= \rll^{\xi_0}\cdot D_\theta(\xi,\eta)\cdot \rrr^{\xi_{N+1}}.
\ee
Hence, the formula for the non-equilibrium stationary state \eqref{corform} becomes:
\be\label{corform1}
\int \widehat D_\theta(\xi, \eta) \mu_\theta^{\text{st}} (d\eta) = \rrr^{|\xi|}\cdot \E_\xi \Bigg[\(\frac{\rll}{\rrr}\)^{\xi_0(\infty)} \Bigg]=  \rrr^{|\xi|}\cdot G\(\xi, \frac{\rll}{\rrr}\)
\ee
where $G(\xi,\cdot)$ is the generating-function defined in Section \ref{TWO}, equation \eqref{genplace}.
The following theorem then follows by combining \eqref{corform1} with Theorem \ref{GEE} and expresses the difference between the expectations of the duality functions in the non-equilibrium stationary measure and their non-interacting counterparts in terms of the probabilities of $\kappa$ dual particles being all absorbed at one end.
\bt \label{Teo:MoM} Let $\bix\in V_n$, with $n\ge 2$, then
\beq\label{MoM}
\int \widehat D_\theta (\phi(\bix), \eta) \mu_\theta^{\text{st}} (d\eta)-\int \widehat D_0 (\phi(\bix), \eta) \mu_0^{\text{st}} (d\eta)=
 \sum_{\kappa=2}^n\gamma_\kappa(\bix)\cdot
 \(\rrr-\rll\)^\kappa \rll^{n-\kappa}
\eeq
with
\beq\label{gamm}
\gamma_\kappa(\bix) &:=& \sum_{\biy\in C_\kappa(\bix)} \gee(\phi(\biy),0)\nonumber\\
\gee(\xi,0)&=&\pee_{\xi}(\xi_0(\infty)=0)-\pee^{\text{irw}}_{\xi}(\xi_0(\infty)=0)
\eeq
and $C_\kappa(\bix)$ as  in \eqref{cmx}.
\et
\noindent

The following corollary of Theorem  \ref{Teo:MoM} specializes to expectations of products of occupation numbers at different sites in the non-equilibrium steady state.
\bc\label{COR}
For every
$\bix= (x_1,\ldots, x_n)\in V_n$ such that $x_i\not= x_j$ for all $i,j\in \{1,\ldots,n\}, i\not= j$,  we have
\beq\label{CoR3}
\E_{\mu^{\text{st}}_\theta}\bigg[\prod_{i=1}^n \eta_{x_i}\bigg]- \prod_{i=1}^n\E_{\mu^{\text{st}}_\theta}\left[\eta_{x_i}\right]= \sum_{\kappa=2}^n\gamma_\kappa(\bix)\cdot
 \(\rrr-\rll\)^\kappa \rll^{n-\kappa}
\eeq
where
$\E_{\mu^{\text{st}}_\theta}\left[\eta_{x}\right]=\sum_{y\in }  \rho_y \, q(x,y)$ with $q(\cdot,\cdot)$ as in \eqref{ABS} and
$\gamma_\kappa(\bix)$ as in \eqref{gamm}.
\ec
\bpr Using \eqref{onedualpart}, \eqref{prod}  and Theorem \ref{Teo:MoM}, we have for all $\theta$
\beq
\int \widehat D_\theta (\phi(\bix), \eta) \mu_\theta^{\text{st}} (d\eta)=\E_{\mu^{\text{st}}_\theta}\bigg[\prod_{i=1}^n \eta_{x_i}\bigg]
\eeq
and then, in particular, for $\theta=0$, from the product-nature of $\mu^{\text{st}}_0$  we have
\beq
\int \widehat D_0 (\phi(\bix), \eta) \mu_0^{\text{st}} (d\eta)=\prod_{i=1}^n \E_{\mu^{\text{st}}_0}\left[\eta_{x_i}\right]
\eeq
This concludes the proof.
\epr
\vskip.2cm
\noindent

\noindent
\br
Specializing the corollary to $n=2$ we obtain information on the covariance. We see that the
explicit dependence on the boundary densities $\rll$ and $\rrr$ that turns out to be a quadratic  function of their difference: for $x\neq y$,
\beq\label{CoV}
\text{cov}_{\mu^{\text{st}}_\theta}(\eta_x,\eta_y)&:=&\E_{\mu^{\text{st}}_\theta}\left[\eta_{x}\eta_y\right]- \E_{\mu^{\text{st}}_\theta}\left[\eta_{x}\right]\cdot \E_{\mu^{\text{st}}_\theta}\left[\eta_{y}\right]\nn\\
&=&   \({\rrr}-\rll\)^2 \cdot \gee(\delta_x+\delta_y,0)\nn\\
&=&  \({\rrr}-\rll\)^2 \cdot \(\pee_{\delta_x+\delta_y}(\xi_0(\infty)=0)-\pee^{\text{irw}}_{\delta_x+\delta_y}(\xi_0(\infty)=0)\).
\eeq
The covariance
is {\em exactly} quadratic in $(\rll-\rrr)$ with a multiplying factor (i.e. the difference of the two absorption probabilities above)  not depending on $\rll$ and $\rrr$.  This multiplying factor is non-positive for exclusion particles, $\theta<0$ (by Liggett's inequality \cite{lig}, chapter 8) and non-negative for inclusion particles, $\theta>0$
(by the analogue of Ligget's inequality from \cite{grv}).
\er

\subsection{Examples}
In this section we consider three concrete examples where the non-equilbrium steady state expectations are known in closed form, and we explicitly  verify the recursion \eqref{dualrecu}.

We restrict to the situation  where the lattice $V=\{1,\ldots, N\}$  is viewed as a one-dimensional chain and particles jump only to nearest neigbors, and interact (in a symmetric way) only if sitting   in neighboring sites. More precisely we chose the function $p(\cdot,\cdot)$ in \eqref{GenRes1} and \eqref{GenDual1} as
\beq
p(i,j)=\mathbf 1_{j=i\pm 1}
\eeq
Moreover we chose the reservoirs clocks to have rates 1:
\beq
\cll=\crr=1
\eeq
\subsubsection*{Independent  walkers (case $\theta=0$)}
In this case the absorption probabilities of   the dual-process:
\beq\label{GenDual2}
\loc_0^{\text{Dual}} f(\xi)&:=& \sum_{i=1}^ N\, \xi_i[f(\xi^{i,i+1})-2f(\xi)+f(\xi^{i,i-1})]\nn
\eeq
 can be explicitly computed and are determined by the single-walker absorption probabilities which are given by:
\beq\label{px}
q^+_x:=q(x,N+1)={\bf P}_x(X^{\text{rw}}(\infty)=N+1)= \frac x{N+1}, \qquad x\in\{0, 1, \ldots, N+1\}.
\eeq
Here ${\bf P}_x$ is the probability law of the random walker $\{X^{\text{rw}}(t),\, t\ge 0\}$ on $V^*=\{0,1,\ldots,N+1\}$, starting from $x\in V^*$, with rate-one nearest-neighbour jumps and with $\{0,N+1\}$ absorbing states.
\vskip.1cm
\noindent
 If we now consider the corresponding system with reservoirs:
 \begin{eqnarray}\label{GenRes2}
{\loc_0^{\text{res}}}f(\eta)&=& \sum_{i=1}^{N-1} \eta_i[f(\eta^{i,i+1})-f(\eta)]+ \sum_{i=2}^{N} \eta_i[f(\eta^{i,i-1})-f(\eta)]\\
&+&\left\{\rll[f(\eta+\delta_1)-f(\eta)]+  \eta_1 [f(\eta-\delta_1)-f(\eta)]\right\}\nonumber\\
&+&\left\{\rrr [f(\eta+\delta_{N})-f(\eta)]+ \eta_{N} [f(\eta-\delta_{N})-f(\eta)]\right\}\nonumber
\end{eqnarray}
 we have that the stationary measure $\mu_0^{\text{st}}$  is a non-homogeneous product measure with Poisson-distributed marginals:
 \beq\label{mu0}
\mu_0^{\text{st}}\sim \otimes_{x=1}^{N} \text{Pois}\(\rll +(\rrr-\rll) q^+_x\)
\eeq
with $q^+_x$ as in \eqref{px}.

\vskip.3cm
\noindent
Notice that in this case, since we can easily  compute the generating function $G^{\text{irw}}(\cdot,\cdot)$, it is possible to directly verify the  recursion relation \eqref{bonan}. For $\xi\in \Omega^*$ we have
\[
G^{\text{irw}}(\xi, z)= \prod_{i=0}^{N+1} \left(z \left(1-q^+_i\right) + q^+_i\right)^{\xi_i}.
\]
As a consequence,
\begin{eqnarray*}
&&(1-z)\, \frac d {dz}G^{\text{irw}}(\xi, z)
 =  (1-z)\sum_{j=0}^{N+1}\bigg\{\xi_j (1-q^+_j) \big(z (1-q^+_j ) + q^+_j\big)^{\xi_j-1}
\prod_{\substack{i=0 \\ i\not= j}}^{N+1} \big(z (1-q^+_i) + q^+_i\big)^{\xi_i}\bigg\}\nn
\\
&&=
\sum_{j=1}^N\bigg\{\xi_j \(-z \(1-q^+_j \) - q^+_j +1\)\(z \(1-q^+_j \) + q^+_j\)^{\xi_j-1}
 \prod_{\substack{i=1\\ i\not= j}}^N \(z (1-q^+_i ) + q^+_i\)^{\xi_i})\bigg\}\\
 &&=
-|\xi|\, G^{\text{irw}}(\xi, z) + \sum_{j=0}^{N+1}  \xi_j G^{\text{irw}}(\xi-\delta_j, z)
\end{eqnarray*}
and then  \eqref{bonan} is satisfied.
Since we explicitly know both the generating function $G^{\text{irw}}(\cdot,\cdot)$ and the stationary measure \eqref{mu0},   it is possible to verify a posteriory  the duality relation \eqref{corform1}. It is also possible to verify that,
by iterating the recursion relation in its integrated form \eqref{dualrecu} one can recover the    generating function $G^{\text{irw}}(\xi,z)$ starting from the  knowledge of $G^{\text{irw}}(\xi,0)$.

\br
As a further application  of the recursion \eqref{bonan} in the same spirit, one can easily show by induction that if the probabilities for all the particles to be absorbed at zero factorize, i.e., if $G(\xi, 0)= \prod_{i} G(\delta_i,0)^{\xi_i}$ for all $\xi$, then the generating function factorizes, i.e., $G(\xi, z)= \prod_{i} G(\delta_i,z)^{\xi_i}$
and as a consequence the system has a product invariant non-equilibrium stationary measure $\mu^{\text{st}}_0$.
\er

\subsubsection*{Interacting walkers   (case $\theta\neq 0$)}

%

For the interacting case,  we consider two special cases as example. The first example is the simple exclusion process. This is the only interacting model in the class for which there is a full knowledge of the $n$-points correlation functions. The second example  is a special case of the inclusion process  for which an exact formula is known only for the two-point correlations \cite{GKR}.

For the interacting case, specialized to nearest neighbor jumps  and rate 1 reservoirs, the process has generator
\begin{eqnarray}\label{GenRes3}
{\loc_\theta^{\text{res}}}f(\eta)&=& \sum_{i=1}^{N-1}\, \eta_i(1+\theta\eta_{i+1})[f(\eta^{i,i+1})-f(\eta)]+\sum_{i=2}^{N}\, \eta_i(1+\theta\eta_{i-1})[f(\eta^{i,i-1})-f(\eta)]\nn\\
&+&\left\{\rll(1+\theta \eta_1)[f(\eta+\delta_1)-f(\eta)]+(1+\theta \rll)  \eta_1 [f(\eta-\delta_1)-f(\eta)]\right\}\nonumber\\
&+&\left\{\rrr (1+\theta \eta_{N})[f(\eta+\delta_{N})-f(\eta)]+(1+\theta \rrr)  \eta_{N} [f(\eta-\delta_{N})-f(\eta)]\right\}\nonumber
\end{eqnarray}

\subsubsection*{Simple exclusion process (case $\theta=-1$)}
For exclusion process the matrix product ansatz gives an algebraic procedure to calculate all
correlation functions. This provides a recursion relation for the correlation functions
(formula (A.7) in \cite{derrida2007entropy}) that reads
\begin{eqnarray}
\label{derrida}
\E^{\scriptscriptstyle{\rm N}}_{\mu^{\text{st}}}[\eta_{x_1}\eta_{x_2}\ldots \eta_{x_m}]
& = &
(\rll-\rrr)\Big( 1 - \frac{x_{m}}{N+1} \Big )\,
\E^{\scriptscriptstyle{\rm N-1}}_{\mu^{\text{st}}}[\eta_{x_1}\eta_{x_2}\ldots \eta_{x_{m-1}}] \nonumber \\
& + &
\rrr \,\E^{\scriptscriptstyle{\rm N}}_{\mu^{\text{st}}}[\eta_{x_1}\eta_{x_2}\ldots \eta_{x_{m-1}}]
\end{eqnarray}
where $\E^{\scriptscriptstyle{\rm N}}_{\mu^{\text{st}}}$  denotes expectation
in the non-equilibrium steady states of a system of size $N$ and $\mu^{\text{st}}=\mu^{\text{st}}_{-1}$.
In this section we fix  $\xi\in \Omega_{m}$, with $\xi=\phi(\bix)$, and $\bix=(x_1,\ldots, x_m)$ with $1\le x_1 <x_2 < \ldots < x_m\le N$
and denote by $q^{\scriptscriptstyle{\rm (N)}}_{\xi}(k)$  the absorption probabilities:
\beq
q^{\scriptscriptstyle{\rm (N)}}_{\xi}(k):=\pee_\xi(\xi_0(\infty)=k), \qquad k\in \{0,1,\ldots, m\}
\eeq
in a system of size $N$.
Duality yields
\begin{eqnarray}
\label{dual-sep}
\E^{\scriptscriptstyle{\rm N}}_{\mu^{\text{st}}}[\eta_{x_1}\eta_{x_2}\ldots \eta_{x_m}]
& = &
\rrr^m \sum_{k=0}^{m} \(\frac{\rll}{\rrr}\)^{k} \, q^{\scriptscriptstyle{\rm (N)}}_{\xi}(k)
\end{eqnarray}
Inserting \eqref{dual-sep} in \eqref{derrida} the principle of identity of polynomials turns
the recurrence relation for the correlation functions into a recurrence relation for the absorption probabilities:
\begin{eqnarray}
\label{recurrence}
q^{\scriptscriptstyle{\rm (N)}}_{\xi}(k) =
{\mathbf 1}_{k\neq 0}\cdot  q_{\xi-\delta_{x_m}}^{\scriptscriptstyle{\rm (N-1)}}(k-1)\cdot q_{\delta_{x_m}}^{\scriptscriptstyle{\rm (N)}}(1)
+
{\mathbf 1}_{k\neq m}\cdot \Big[q_{\xi-\delta_{x_m}}^{\scriptscriptstyle{\rm (N)}}(k) - q_{\xi-\delta_{x_m}}^{\scriptscriptstyle{\rm (N-1)}}(k)\cdot q_{\delta_{x_m}}^{\scriptscriptstyle{\rm (N)}}(1) \Big]
\end{eqnarray}
where $k \in \{0,1,2,\ldots, m\}$.
Introducing the generating function
\be
G^{\scriptscriptstyle{\rm (N)}}(\xi,z) = \sum_{k=0}^{m}  z^{k} \,q_{\xi}^{\scriptscriptstyle{\rm (N)}}(k)
\ee
the recursion relation of the absorption probabilities \eqref{recurrence} implies
the recursion relation for the generating function
\be
\label{recurr-g}
G^{\scriptscriptstyle{\rm (N)}}(\xi,z) = (z-1)\;  q_{\delta_{x_m}}^{\scriptscriptstyle{\rm (N)}}(1) \cdot G^{\scriptscriptstyle{\rm (N-1)}}(\xi-\phi(x_m),z) + G^{\scriptscriptstyle{\rm (N)}}(\xi-\phi(x_m),z)
\ee
Clearly, if $m=1$, namely $\xi=\delta_x=\phi(x)$ for some  $x\in V$ we have
\be
\label{g1}
G^{\scriptscriptstyle{\rm (N)}}(\phi(x),z) = \frac{x}{N+1} + \Big(1- \frac{x}{N+1}\Big)z
\ee
Thus for the exclusion process the probability generating function for the number
of particles absorbed at zero can be computed by iterating \eqref{recurr-g}.

\noindent
For $m=2$ and $x<y$ the recurrence \eqref{recurr-g} gives
\begin{eqnarray}
G^{\scriptscriptstyle{\rm (N)}} (\phi(x,y),z)
=
(z-1) \; q_{\delta_{y}}^{\scriptscriptstyle{\rm (N)}}(1) \cdot G^{\scriptscriptstyle{\rm (N-1)}}(\phi(x),z) + G^{\scriptscriptstyle{\rm (N)}}(\phi(x),z)
\end{eqnarray}
Using \eqref{g1} we get
\begin{eqnarray}
\label{g2}
G^{\scriptscriptstyle{\rm (N)}} (\phi(x,y), z)
=
(z-1)\Big(1- \frac{y}{N+1}\Big)\Big[\frac{x}{N} +  \Big(1- \frac{x}{N}\Big)z \Big]
+
\frac{x}{N+1} + \Big(1- \frac{x}{N+1}\Big)z
\end{eqnarray}
and one can check that this expression satisfies \eqref{bonan} with $m=2$, i.e.
\be
(1-z)\frac{d}{dz} G^{\scriptscriptstyle{\rm (N)}}(\phi(x,y), z)  + 2 G^{\scriptscriptstyle{\rm (N)}}(\phi(x,y), z) = G^{\scriptscriptstyle{\rm (N)}}(\phi(x), z) + G^{\scriptscriptstyle{\rm (N)}}(\phi(y), z)
\ee

\noindent
For $m=3$ and $x<y<u$ the recurrence  \eqref{recurr-g} gives
\begin{eqnarray}
G^{\scriptscriptstyle{\rm (N)}} (\phi(x,y,u),z)
=
(z-1) \; q_{\delta_{u}}^{\scriptscriptstyle{\rm (N)}}(1) \cdot G^{\scriptscriptstyle{\rm (N-1)}}(\phi(x,y),z) + G^{\scriptscriptstyle{\rm (N)}}(\phi(x,y),z)\nn
\end{eqnarray}
Using \eqref{g2} we get
\begin{eqnarray}
G^{\scriptscriptstyle{\rm (N)} }(\phi(x,y,u),z)
&=&
(z-1)^2\Big(1- \frac{u}{N+1}\Big)
\Big(1- \frac{y}{N}\Big)\Big[\frac{x}{N-1} +  \Big(1- \frac{x}{N-1}\Big)z \Big]
\nonumber\\
&+&
(z-1)\Big(2- \frac{u+y}{N+1}\Big)
\Big[
\frac{x}{N} + \Big(1- \frac{x}{N}\Big)z
\Big] \nonumber\\
&+ &
\frac{x}{N+1} + \Big(1- \frac{x}{N+1}\Big)z
\end{eqnarray}
One can check that this expression satisfies \eqref{bonan} with $m=3$, i.e.
\beq
&&(1-z)\;\frac{d}{dz} G^{\scriptscriptstyle{\rm (N)}}(\phi(x,y,u),z)  + 3 G^{\scriptscriptstyle{\rm (N)}}(\phi(x,y,u),z) \nn\\
&&\hskip4cm= G^{\scriptscriptstyle{\rm (N)}}(\phi(x,y),z) + G^{\scriptscriptstyle{\rm (N)}}(\phi(y,u),z) + G^{\scriptscriptstyle{\rm (N)}}(\phi(u,y),z)
\eeq

\subsubsection*{Inclusion process with $\theta=2$}
For the inclusion process SIP$(2)$ we can verify \eqref{bonan} only for the case  $m=2$ using the results
in Section 5 of  \cite{GKR}. Writing out
\be
G(\phi(x,y),z) = q_{\phi(x,y)}(0) + z q_{\phi(x,y)}(1)  + z^2 q_{\phi(x,y)}(2)
\ee
we see that \eqref{bonan} is equivalent to
\be
\label{verify}
2 q_{\phi(x,y)}(0) +  q_{\phi(x,y)}(1) = \frac{x+y}{N+1}
\ee
From  Eq. (5.6) in \cite{GKR} giving the two-point correlation function for the Brownian
Energy process with reservoirs, and using the fact that the inclusion process with absorbing boundaries (see eq. (3.2) in \cite{GKR}) is dual to it, we can read off
the absorption probabilities as:
\be
q_{\phi(x,y)}(0) = \frac{x(2+y)}{(N+1)(N+3)}
\ee
and
\be
q_{\phi(x,y)}(1) = 1 - \Big(1-\frac{x}{N+3}\Big)\Big(1-\frac{y}{N+1}\Big) - \frac{x(2+y)}{(N+1)(N+3)}
\ee
Thus equation \eqref{verify} is verified.
\\[1cm]
{\large{\bf Acknowledgement:}} We thank Mario Ayala and Federico Sau for useful discussions and suggestions.

\end{document}